\documentclass[11pt]{article}

\usepackage{lineno, hyperref}
\modulolinenumbers[1]

\usepackage{color}
\usepackage{latexsym}
\usepackage{dsfont}
\usepackage{amssymb}
\usepackage{float}
\usepackage{algorithm,algorithmic}
\usepackage{graphicx}
\usepackage{amsmath, amsfonts,amssymb,theorem,euscript,array,enumerate,amsfonts,mathrsfs}
\usepackage{subcaption}
\usepackage{ulem}

\definecolor{newroyal}{rgb}{0.6,0.4,0.8}

\newtheorem{Theorem}{Theorem}[section]
\newtheorem{Definition}[Theorem]{Definition}
\newtheorem{Proposition}[Theorem]{Proposition}

\newtheorem{Lemma}[Theorem]{Lemma}
\newtheorem{Corollary}[Theorem]{Corollary}
\newtheorem{Remark}[Theorem]{Remark}

\numberwithin{equation}{section}

\theorembodyfont{\upshape}
\newcommand{\nc}{\newcommand}
\nc{\ind}{\mathds{1}}
\def \trans{^{\scriptscriptstyle{\intercal}}}

\newcommand{\R}{\mathbb{R}}

\newcommand{\E}{\mathcal{E}}
\newcommand{\F}{\mathcal{F}}

\DeclareMathOperator{\esssup}{esssup}

\def\esssup_#1{\underset{#1}{\mathrm{ess\,sup\, }}}
\def\essinf_#1{\underset{#1}{\mathrm{ess\,inf\, }}}
\def\argmax_#1{\underset{#1}{\mathrm{arg\,max\, }}}
\def\argmin_#1{\underset{#1}{\mathrm{arg\,min\, }}}

\def\reff#1{{\rm(\ref{#1})}}
\def \ep{\hbox{ }\hfill$\Box$}

\def \Inf{\displaystyle\inf}

\def \Max{\displaystyle\max}

\def\b1{\bf 1}

\def \N{\mathbb{N}}
\def \R{\mathbb{R}}

\def \E{\mathbb{E}}
\def \F{\mathbb{F}}

\def \L{\mathbb{L}}
\def \P{\mathbb{P}}

\def \Ac{{\cal A}}

\def \Cc{{\cal C}}
\def \Dc{{\cal D}}

\def \Fc{{\cal F}}
\def \Gc{{\cal G}}
\def \Hc{{\cal H}}

\def \Kc{{\cal K}}

\def \Pc{{\cal P}}
\def \Mc{{\cal M}}

\def \Tc{{\cal T}}
\def \Uc{{\cal U}}

\def \ep{\hbox{ }\hfill$\Box$}

\def \trans{^{\scriptscriptstyle{\intercal}}}

\def\reff#1{{\rm(\ref{#1})}}
\def\beqs{\begin{eqnarray*}}
\def\enqs{\end{eqnarray*}}
\def\beq{\begin{eqnarray}}
\def\enq{\end{eqnarray}}

\begin{document}



\title{It\^o's formula for flows of measures on semimartingales\thanks{We are grateful to the authors of \cite{taltouzha21}, who  pointed out in private communication  an error  in the first version of this paper, and obtained in their  work a similar It\^o formula by different techniques.}}

\author{Xin GUO\footnote{Department of Industrial Engineering and Operations Research, University of California, Berkeley, USA \sf xinguo at berkeley.edu} 
\quad\quad
Huy\^en PHAM\footnote{
LPSM, Universit\'e  Paris Cit\'e,  Building Sophie Germain, Avenue de France, 75013 Paris, \sf pham at lpsm.paris  
The work of this author  is supported by FiME (Finance for Energy Market Research Centre), and the ANR-18-IDEX-0001.
}
\quad\quad 
Xiaoli WEI\footnote{Tsinghua-Berkeley Shenzhen Institute, Tsinghua University, 518055 Shenzhen \sf xiaoli\_wei@sz.tsinghua.edu.cn}}

 \date{First version: October 13, 2020 \\
This version: September 17, 2022}

\maketitle
 
\begin{abstract}
We establish It\^o's formula along  flows of probability measures associated with general semimartingales; this generalizes existing results for flows of measures on It\^o processes. Our approach is to first establish It\^o's formula for cylindrical functions and then extend it to the general case via function approximation and localization techniques.

This general form of It\^o's formula enables the  derivation of  dynamic programming equations and  verification theorems
 for  McKean--Vlasov controls with jump diffusions and for McKean--Vlasov mixed regular-singular control problems.
 It also allows for generalizing the classical relationship between the maximum principle and the dynamic programming principle to the  McKean--Vlasov singular control setting, where the adjoint process  is expressed in terms of the derivative of the value function with respect to the probability measures.

\end{abstract}




\vspace{5mm}


\vspace{5mm}







\section{Introduction}


It\^o's formula is one of the most fundamental building blocks in stochastic calculus. It is the key tool to connect PDEs with functional expectation of diffusion processes.
Recent developments in mean-field games and mean-field controls have extended It\^o's formula to flows of probability measures.  Analogously, this extension is a natural tool for deriving PDEs in the space of probability measures when dealing with SDEs of mean-field type called McKean--Vlasov equations.
For instance, the generalized  It\^o's formula has been used to derive the  master equations in mean-field games and the Bellman dynamic programming equation for McKean--Vlasov control problems. (See \cite{CD2018I} and the references therein.)

There are three approaches to establish It\^o's formula along flows of measures associated with It\^o processes. The first is to discretize  time and mimic the standard proof of It\^o's formula~\cite{BLPR2014}; this approach can also yield the  It\^o's formula for flows of measures associated with mean-field jump diffusions~\cite{Li2018}.
The second is to use the Fokker--Planck approach by assuming the existence of the density for the probability measure~\cite{CDLL2019}.
 The third, known as the particle approximation approach, is to approximate  flows of measures by flows of empirical measures~\cite{CDD2014,  CD2018I}; this has been applied recently 
 in ~\cite{RP2019}, \cite{RP2021} to derive several It\^o--Wentzell--Lions formulae on Wiener spaces for real-valued functional random fields that depend on measure flows.

\paragraph{McKean--Vlasov jump diffusion and singular control} McKean--Vlasov processes, first introduced and studied in \cite{MCKEAN1969}, are stochastic processes governed by SDEs whose coefficients depend on distributions of the solutions. McKean--Vlasov controls are concerned with the optimal control of McKean--Vlasov processes.  As in classical control theory, the two main tools to analyze McKean--Vlasov controls are the stochastic maximum principle and the dynamic programming principle.
 The former has been developed in~\cite{anddje10, CD2015} to study  controlled McKean--Vlasov systems in terms of adjoint backward SDEs; it has also been adopted for mean-field games with singular controls~\cite{FH2017}. The latter has been developed through a series of works  including for Markovian controls~\cite{PW2018}, open-loop controls~\cite{BCP2018, cospha19}, Markovian and non-Markovian frameworks~\cite{DPT2019}, and closed-loop controls~\cite{WZ2020}.  None of these works except for ~\cite{FH2017}, however, involves jumps.

The dynamic programming principle was developed to analyze the viscosity solution property for a class of  jump diffusion  processes whose coefficients and control are independent of  the state process~\cite{BIVS2019}.
 Moreover, most existing works on McKean--Vlasov controls deal with It\^o processes  with continuous paths, except for~\cite{HOS2017, FH2017, HSSD2018}, which considered singular controls and~\cite{HAS2014, BIVS2019}, which studied jump diffusions.
 To the best of our knowledge, there is no prior work on
  the dynamic programming approach for McKean--Vlasov controls with general semimartingales.  The  barrier has been the lack of a general form of It\^o's formula  for flows of measures on semimartingales.

\paragraph{Our work}  We establish It\^o's formula for flows of measures associated with general and possibly discontinuous semimartingales (Theorems~\ref{TheoremIto} \& \ref{Theoremito_diffcond}).
It enables us to derive dynamic programming equations and  verification theorems for  McKean--Vlasov controls with jump diffusions and for McKean--Vlasov mixed regular-singular control problems (Theorems \ref{Propjumpver} \& \ref{veri}).
 It also allows for generalizing the classical relationship between the maximum principle and the dynamic programming principle  to the  McKean--Vlasov mixed regular-singular control setting, where the adjoint process is expressed in terms of the derivative of the value function with respect to probability measures along the optimal path (Theorem \ref{Thmrelationadjoint}).

Our approach starts by establishing  It\^o's formula on  cylindrical functions  \cite{Pro05}. Then, by invoking a general form of Stone--Weierstrass theorem on compact sets in the joint Wasserstein space and Euclidean space,  this class of mean-field functions is shown to be dense in the class of twice differentiable  functions on the space of probability measures of order two with the topology of $\Cc^{1, 1}$.  Finally, with an appropriate localization argument, the general form of It\^o's formula is  established.

 There are several key technical ingredients in our work. The first are the properties of cylindrical functions. These  are smooth mean-field functions with integrable forms,  initially studied in Fleming--Viot processes~\cite{FV1979} for modeling population genetics  and further developed for analyzing a general class of probability measure-valued processes called polynomial diffusions~\cite{CLSF2019}.
The second
are linear derivatives on the space of probability measures~\cite{FV1979, CLSF2019};
this form of derivatives enables characterizing the behavior of both the jumps of flows of probability measures and the jumps of  semimartingale processes. This is in contrast to  the  It\^o's formula for flows of measures on diffusion processes,  which involves only  Lions derivatives~\cite{Lions2012}.  Indeed, Lions derivatives, obtained by identifying the Wasserstein space with the Hilbert space of random variables on atomless probability space~\cite{BLPR2014, CDD2014},  appear insufficient for characterizing the infinitesimal changes on the functional of controlled McKean--Vlasov  processes with the addition of jumps; this has been observed before in~\cite{BIVS2019}.

\paragraph{Comparison with related works} In an independent work, which has been brought to our attention during the review of this paper, ~\cite{taltouzha21}
obtained a similar form of It\^o's formula for flows of measures on c\`adl\`ag semimartingales.  There are two main differences  between our work and theirs. First,  they use a time-discretization approach, which is different from our cylindrical function technique.  Secondly, the conditions for the It\^o's formula in their work are different from ours. They are essentially different sets of conditions needed for the dominated convergence theorem to allow for the interchange of the expectation and integration; see Remark~\ref{comaretaltouzha21} for more detailed comparison.

\paragraph{Cylindrical approach} The cylindrical function technique  appears natural for handling discontinuous semimartingales and is of independent interest.
It has recently been used
to obtain It\^o's formula for controlled measure-valued martingales~\cite[
Theorem 5.1]{CKLS2021}.  Cylindrical functions have also been used to study viscosity solutions for controlled McKean--Vlasov dynamics~\cite{BIVS2019}.

Besides the cylindrical function technique and the time-discretization approach, it is naturally of interest to explore for instance the particle approximation method. However, it appears  more difficult through this route for the case of discontinuous semimartingales, and we leave it for future studies.
(See Remark~\ref{remark:particle} for more technical discussions on the  particle approximation attempt).

\paragraph{Outline of the paper} Section \ref{secnot} introduces the notation  and recalls definitions of several forms of derivatives in the Wasserstein space of probability measures.  Section~\ref{secIto}  is devoted to  It\^o's formula and  its variants for flows of probability measures for semimartingales.
Section~\ref{secMKVjump} studies McKean--Vlasov control with jump diffusions and
Section~\ref{secMKVsingular} analyzes McKean--Vlasov mixed regular-singular controls.

\section{Notations and Preliminaries} \label{secnot}

Throughout the paper,  $(\Omega, \Fc, \F=(\Fc_t)_{t \geq 0}, \P)$
is a fixed filtered probability space satisfying the usual conditions. That is,  $(\Omega, \Fc, \F, \P)$ is complete and the filtration is right continuous.
We assume that there exists a sub-$\sigma$-algebra $\Gc$ of $\Fc$, with  $\Fc$ rich enough as will be explained shortly.
Moreover, we will adopt the following notations, unless otherwise specified.
\begin{itemize}
\item  Given any normed space $(E,|\cdot|)$,  $\mathcal{P}  (E) $ is the set of all probability measures on $(E,|\cdot|)$. For any $p \in \mathbb{N}$, $\mathcal{P}_p (E) $ is the set of all probability measures of $p$-th order on $(E,|\cdot|)$, defined as   $\mathcal{P}_p (E) = \biggr\lbrace \mu \in \mathcal{P}  (E) \biggr|\|\mu\|_p: = \biggr(\int_E |x|^p  \mu( dx)\biggr)^\frac{1}{p}  < \infty \biggr\rbrace.$ For instance,
 $\mathcal{P}_2(E) $ is the set of all probability measures with finite second-order moment.  For any  probability measures $\mu, \mu'$   in  $\mathcal{P}_p(E)$, the $p$-th order Wasserstein distance on $\mathcal{P}_p(E)$ is  defined as
$W_p(\mu,\mu')=  \inf_{\pi}   \limits  \left(\int_{E\times E}  |y-y'|^p  {\pi} (dy,dy')  \right)^{\frac{1}{p}},$
where $\pi $ is a coupling of $\mu $ and $\mu' $ in the sense  that $\pi \in \mathcal{P}  (E\times E)$ with marginals $\mu$ and $\mu'$. For any probability measure $\mu \in \Pc(E)$, ${\rm supp}(\mu)$ is the support of $\mu$.  $L^2(E)$ is the space of all square-integrable functions $f: E \to \R$. 

\item 
$L^2_\mu(\R^d; \R^d)$ is the space of all Borel mappings $f: \R^d \to \R^d$ such that $\int_{\R^d} |f|^2\mu(dx) < \infty$. $\nabla \Cc_c^\infty (\R^d):=\{\nabla \varphi: \varphi \in \Cc^\infty_c(\R^d)\}$, where $\Cc^\infty_c(\R^d)$ is the space of smooth real functions with compact support in $\R^d$. The closure of $\nabla \Cc_c^\infty (\R^d)$ in $L^2_\mu(\R^d; \R^d)$ is denoted as $\Tc_\mu \Pc_2(\R^d)$, which is called tangent space of $\Pc_2(\R^d)$ at a given measure $\mu \in \Pc_2(\R^d)$.

\item For vectors $a$, $b$ $\in$ $\R^d$, $d \in \N^{+}$, $a.b$ $=$ $\sum_{i=1}^d a_i b_i$. For the matrix $C=(C_{ij})_{1 \leq i, j\leq d}$ $\in$ $\R^{d \times d}$, ${\rm Tr}(C)$ $=$ $\sum_{i=1}^d C_{ii}$, the transpose of $C$ is $C\trans$. Denote ${\rm diag}(a_1, \ldots, a_d)$ for the diagonal matrix in $\R^{d \times d}$ with diagonal entries $a_i \in \R$, $1 \leq i \leq d$.


\end{itemize}
We denote  $L^2(\mathcal{F}; E$) as the space of all $E$-valued square integrable random variables on $(\Omega, \Fc, \P)$; for any $\vartheta \in L^2(\Fc; E)$, we set $\|\vartheta\|_{L^2} = \E[|\vartheta|^2]^{\frac{1}{2}}$. We will assume that $\Gc$ is "rich enough" in a sense that $\Pc_2(E) = \{\P_{\xi}; \xi \in L^2(\Gc; E)\}$, where
$\P_{\xi}$ denotes the law of $\xi$. This is satisfied whenever the probability space
$(\Omega, \Gc, \P)$ is atomless, see page 352, \cite{CD2018I}.
We may use $(\bar{\Omega}, \bar{\mathcal{F}},\bar{\P})$ for an independent and identical copy of the probability space $(\Omega, \mathcal{F},  \P)$, and $\bar{\vartheta} \in  L^2(\bar{\mathcal{F}}; {E})$ for an independent and identical copy of a random variable $\vartheta \in  L^2( \mathcal{F};  E)$.

\paragraph{Differentiability of functions of probability measures}

We first recall two notions of differentiability of functionals with respect to probability measures that will be used throughout the paper. One is the Lions derivative from the identification of the Wasserstein space with the Hilbert space of random variables on an atomless probability space (\cite{Lions2012} and \cite{CD2018I}), the other is the standard/ linear derivative on the Wasserstein space  (\cite{FV1979}, \cite{CLSF2019}, and \cite{CD2018I}).

The Lions derivative of a functional $f$, introduced in \cite{Lions2012}, is defined through the lift of $f$. 
The idea is to take the function $f:\mathcal{P}_2(\mathbb{R}^d ) \rightarrow \mathbb{R} $, and let $\tilde{f}:L^2( \mathcal{F} ;  \mathbb{R}^d ) \rightarrow \mathbb{R} $  be a lift  of $f$ such that $\tilde{f} ( \vartheta) = f( \P_\vartheta) $ for any $\vartheta \in L^2( \mathcal{F} ;  \R^d)$. Then, $\tilde{f}$ is differentiable in the Fr\'echet sense at $\vartheta_0$ if there exists a linear continuous mapping $D\tilde{f} (\vartheta_0):L^2( \mathcal{F} ;  \mathbb{R}^d ) \rightarrow \mathbb{R}$ such that
\begin{align*}
\tilde{f}(\vartheta) - \tilde{f}(\vartheta_0) = \E \big[ D\tilde{f} (\vartheta_0).(\vartheta -\vartheta_0)\big] + o(\|\vartheta -\vartheta_0\|_{L^2}),
\end{align*}
as  $\|\vartheta -\vartheta_0\|_{L^2} \rightarrow 0$.
It has been shown that when $\tilde f$ is the lift of a function $f$ in $\Pc_2(\R^d)$, the law of  $D\tilde f(\vartheta_0)$ depends on $\vartheta_0$ only via its law $\P_{\vartheta_0}$, and
\beq \label{defderL}
D\tilde f(\vartheta_0) &=& h_0(\vartheta_0),
\enq
for some Borel function $h_0$ $:$ $\R^d$ $\rightarrow$ $\R^d$. (See e.g., \cite{CD2018I}, Chapter 5 and \cite{GT2019}). The Lions derivative is thus well defined:  
\begin{Definition} \label{Lions}
$f$ is differentiable at $\mu_0$ $=$ $\P_{\vartheta_0} \in \mathcal{P}_2(\mathbb{R}^d)$ if its lift function $\tilde f$  is Fr\'echet differentiable at $\vartheta_0$; and in this case, the function $h_0$ in \eqref{defderL} is called the Lions derivative of $f$ at $\mu_0$, and denoted as
$(\partial_\mu f)(\mu_0, .)$.
\end{Definition}


\medskip

The definition of the linear derivative on the Wasserstein space $\Pc_2(\R^d)$ (\cite{FV1979}, \cite{CLSF2019} and \cite{CD2018I}, Chapter 5) is more straightforward:
\begin{Definition} \label{linear}For $f:$ $\Pc_2(\R^d)$ $\to$ $\R$, the linear derivative of $f$ at $\mu$ is a function $\frac{\delta f}{\delta \mu}$ on $\Pc_2(\R^d)$ $\times$ $\R^d$ such that for every $\mu$, $\mu'$ $\in$ $\Pc_2(\R^d)$
\beqs
f(\mu) - f(\mu') = \int_0^1 \int_{\R^d} \frac{\delta f}{\delta \mu}(h \mu + (1 - h)\mu', x) (\mu - \mu')(dx) dh.
\enqs
\end{Definition}

Clearly, the linear derivative $\frac{\delta f}{\delta \mu}$ is defined up to an additive constant. It is very closely related to
 the G\^ateaux derivative, denoted by $\frac{\partial f}{\partial\mu}$, which usually assumes the existence of the density function for a probability measure (see \cite{BJS2017}).

Under suitable regularity conditions on $f$, the Lions derivative in Definition \ref{Lions}, and the linear derivative and the G\^ateaux derivative in Definition \ref{linear} are connected through the following equation, according to  Proposition 5.48 of \cite{CD2018I},
\beq \label{reldermu}
\partial_\mu f(\mu, x) = \partial_x \frac{\delta f}{\delta \mu}(\mu, x) = \partial_x \frac{\partial f}{\partial\mu}(\mu, x).
\enq
For instance, consider the linear  $\Phi(\mu)$ $=$ $\int_{\R^d} g(x)\mu(dx)$  $=:$ $\left<g, \mu\right>$, where the derivative of $g$ has linear growth, then  $\partial_\mu \Phi(\mu, x)$ $=$
$\partial_{x} g(x)$ $=$ $\partial_x \frac{\delta  \Phi}{\delta \mu}(\mu,x)$.

\medskip

\begin{Remark}
In addition to the Lions derivative and the linear derivative, there is an intrinsic notion of derivative on the Wasserstein space used in optimal transport and gradient flows (\cite{AGS2008}, \cite{GT2019}). Precisely, one says that  $f: \Pc_2(\R^d) \to \R$ is differentiable at $\mu \in \Pc_2(\R^d)$ if there exists a unique $\xi \in \Tc_\mu \Pc_2(\R^d)$, such that
\beqs
\lim_{n \to \infty} \frac{f(\mu_n) - f(\mu) - \int_{\Pc_2(\R^d) \times \Pc_2(\R^d)} \xi(x)(y-x)\pi_n(dx, dy) }{W_2(\mu, \mu_n)} = 0
\enqs
for any sequence $\{\mu_n\}_n \subset \Pc_2(\R^d)$ such that $2$-Wasserstein distance $W_2(\mu_n, \mu) \to 0$ and any sequence of optimal plans $\{\pi_n\}_n$ $\subset \Pc_2(\R^d \times \R^d)$ with marginals $\mu_n, \mu$ such that $W_2(\mu_n, \mu) = \Big(\int_{\R^d \times \R^d} |x-y|^2 \pi_n(dx, dy)\Big)^{\frac{1}{2}}\}$ for any $n$. Such $\xi$ is called Wasserstein gradient and denoted as $\nabla_w f(\mu)$. It is shown in \cite{GT2019}  that Wasserstein gradient $\nabla_w f(\mu)$ is equivalent to Lions derivative $\partial_\mu f(\mu, \cdot)$ in Definition \ref{Lions}.  
\end{Remark}

Our It\^o's formula (Theorem \ref{TheoremIto}) will be first  established on  the space of differential functions $\mathcal{C}^{1,1} (\mathcal{P}_2(\mathbb{R}^d))$ in $\mathcal{P}_2(\mathbb{R}^d)$, first introduced in \cite{BLPR2014}, \cite{CDD2014}.

\begin{Definition} \label{defregularityf}
We say a  function $f \in \mathcal{C}^{1,1}(\mathcal{P}_2(\mathbb{R}^d))$, if its lift $\tilde f$ is Fr\'echet differentiable, and if there exists a continuous version of $\partial_\mu f(\mu, x)$ such that
\begin{itemize}
\item the mapping $(\mu, x)$ $\mapsto$ $\partial_\mu f(\mu, x)$ is jointly continuous with respect to $(\mu, x)$ and there is a constant $C >0$ such that
    \beq \label{lipsct1}
   |\partial_\mu f(\mu, x)| \leq C,
    \enq
     for any $\mu$ $\in$ $\Pc_2(\R^d)$ and any $x$ $\in$ $\R^d$;
\item for any $\mu$ $\in$ $\Pc_2(\R^d)$, the mapping $x$ $\mapsto$ $\partial_\mu f(\mu, x)$ is continuously differentiable. Its derivative, denoted by $\partial_x\partial_\mu f(\mu, x)$, is jointly continuous with respect to $(\mu, x)$; and there is a constant $C >0$ such that
    \beq \label{lipsct2}
    |\partial_x\partial_\mu f(\mu, x)| \leq C,
    \enq
for any $\mu$ $\in$ $\Pc_2(\R^d)$ and any $x$ $\in$ $\R^d$.
\end{itemize}
\end{Definition}

We also introduce an alternative class
$\Cc^{1, 1}_{relax}(\Pc_2(\R^d))$  of differential functions where we relax   the boundedness conditions in
$\Cc^{1, 1}_{}(\Pc_2(\R^d))$ by square integrability growth conditions, in view of
 It\^o's formula in Theorem \ref{Theoremito_diffcond}.


\begin{Definition}  \label{defregularityf_relax}
We say a  function $f \in \mathcal{C}^{1,1}_{relax}(\mathcal{P}_2(\mathbb{R}^d))$, if its lift $\tilde f$ is Fr\'echet differentiable, and if there exists a continuous version of $\partial_\mu f(\mu, x)$ such that
\begin{itemize}
\item the mapping $(\mu, x)$ $\mapsto$ $\partial_\mu f(\mu, x)$ is jointly continuous with respect to $(\mu, x)$;

\item for any $\mu$ $\in$ $\Pc_2(\R^d)$, the mapping $x$ $\mapsto$ $\partial_\mu f(\mu, x)$ is continuously differentiable. Its derivative, denoted by $\partial_x\partial_\mu f(\mu, x)$, is jointly continuous with respect to $(\mu, x)$;
\item for  any compact set  $\Kc$ in $\Pc_2(\R^d)$,
\beq \label{conditionPhi}
\sup_{\mu \in \Kc} \biggl[\int_{\R^d}\Big \lbrace  |\partial_\mu f(\mu, x)|^2  +  |\partial_x \partial_\mu f(\mu, x)|^2\Big\rbrace\mu(dx) \biggl] < \infty.
\enq
\end{itemize}
\end{Definition}

Note that in the above Definitions, the continuity in $\mu$ is with respect to the Wasserstein distance $W_2$.  Moreover, if $f \in \Cc^{1, 1}(\Pc_2(\R^d))$ or $\Cc^{1, 1}_{relax}(\Pc_2(\R^d))$, then the linear derivative $\frac{\delta f}{\delta \mu}$ in Definition \ref{linear} exists and is jointly continuous with respect to $(\mu, x) \in \Pc_2(\R^d) \times \R^d$.

\section{It\^o's Formula for Flows of Measures on Semimartingales} \label{secIto}

\subsection{It\^o's Formula and Several Variants}

In this section, we will present   It\^o's formula and several of its variants for flows of measures on a class of semimartingales $X= (X_t)_{0 \leq t \leq T}$  ($T > 0$),
which are square integrable and satisfy the {\bf (H)} condition, following the notation in \cite{Pro05}. This condition is similar to the integrability conditions on the drift and diffusion coefficients imposed
for It\^o's diffusion processes in \cite{BLPR2014} and \cite{CDD2014}.

\medskip

\noindent  {\bf (H)} \; There exists one decomposition $X$ $=$ $X(0)$ $+$ $V$ $+$ $L$, where $V$ is c\`adl\`ag process of finite variation with $V(0)$ $=$ $0$, and $L$ is  a local martingale with $L(0)$ $=$ $0$ such that
\beq \label{Xassum}
\E[{\rm Var}(V)_T] < \infty, \;\;\; \E\Big[\sum_{0 < t \leq T} \big|\Delta X_t\big|\Big] < \infty, \;\;\; \E[[X, X]_T^c] < \infty.
\enq
Here ${\rm Var}(V)_T$ denotes the variation of the process $V$ on $[0, T]$,  $X_{t-}$  the left limit of $X$ at $t$, $\Delta X_t = X_t - X_{t-}$  the jump of $X$ at $t$,  $[X, X]$  the quadratic variation of $X$, and $[X, X]^c_t$ $=$ $[X^c, X^c]$ $=$ $[X, X]_t$ $-$ $\sum_{0 \leq r \leq t} |\Delta X_r|^2$  the continuous part of $[X, X]$, with  $X^c$  the continuous part of $X$.


\begin{Theorem}[It\^o's formula] \label{TheoremIto}
Given a semimartingale $X$ satisfying {\bf (H)} and $\Phi : \mathcal{P}_2(\mathbb{R}^d) \rightarrow \mathbb{R}$.
If $\Phi$ $\in$ $\Cc^{1, 1}(\Pc_2(\R^d))$, then for any $0\le t < s \le T $,
\beq  \label{itoPP}
& & \Phi(\P_{X_s }) - \Phi(\P_{X_t})\\
& & =
  \E \biggl[ \int_t^s \partial_\mu \Phi (\P_{X_{r-}}, X_{r-}).d X_r + \frac{1}{2}{\rm Tr}\Big(\partial_x \partial_\mu \Phi(\P_{X_{r-}}, X_{r-}) d[ X,  X]_r^c \Big) \nonumber\\
&  & \; + \sum_{t < r \leq s} \Big\lbrace\Big(\Phi(\P_{X_r}) - \Phi(\P_{X_{r-}})\Big)  1_{\{\P_{X_r} \neq \P_{X_{r-}}\}} - \partial_\mu\Phi(\P_{X_{r-}},  X_{r-}).\Delta X_r  \nonumber\\
&  & \;+ \Big(\frac{\delta \Phi}{\delta \mu}(\P_{X_{r}}, X_{r}) - \frac{\delta \Phi}{\delta \mu} (\P_{X_{r}}, X_{r-}) \Big)  1_{\{\P_{X_r} = \P_{X_{r-}}\}} \Big\rbrace\biggl] \nonumber.
\enq
\end{Theorem}

\begin{Remark} \label{Itoremainder}
Note that when $\Phi$ $\in$ $\Cc^{1, 1}(\Pc_2(\R^d))$, the jump remainder in the RHS of \reff{itoPP} is finite. That is,
\beqs
& & \sum_{t < r \leq s} \Big|\Big(\Phi(\P_{X_r}) - \Phi(\P_{X_{r-}})\Big)1_{\{\P_{X_r} \neq \P_{X_{r-}}\}} \Big| < \infty, \; \E\Big[ \Big|\sum_{t < r \leq s}\partial_\mu\Phi(\P_{X_{r-}},  X_{r-}).\Delta X_r \Big|\Big] < \infty,\\
& &\E\Big[ \sum_{t < r \leq s} \Big|\Big(\frac{\delta \Phi}{\delta \mu}(\P_{X_{r}}, X_{r}) - \frac{\delta \Phi}{\delta \mu} (\P_{X_{r}}, X_{r-})\Big) 1_{\{\P_{X_r} = \P_{X_{r-}}\}}\Big|\Big] <  \infty.
\enqs
Indeed, by the definition of the Lions derivative in Definition \ref{Lions} and its relation with the linear derivative \reff{reldermu}, we see
\beq\label{Phiintegralform}
\Phi(\P_{X_r}) - \Phi(\P_{X_{r-}}) &=& \int_0^1 \frac{d}{d h }\Phi(\P_{X_{r-} + h\Delta X_r}) d h \nonumber\\
&=& \int_0^1 \E\big[\partial_\mu \Phi(\P_{X_{r-} + h \Delta X_r}, X_{r-} + h \Delta X_r). \Delta X_r\big] d h,\\
\label{equ: Itolinearlions}
\frac{\delta \Phi}{\delta \mu}(\P_{X_{r}}, X_{r}) - \frac{\delta \Phi}{\delta \mu} (\P_{X_{r}}, X_{r-}) &=& \int_0^1 \partial_\mu \Phi(\P_{X_r}, X_{r-} + h \Delta X_r).\Delta X_r dh.
\enq
Thus,  by \reff{Phiintegralform}-\reff{equ: Itolinearlions}(or the Lipschitz continuity of $\Phi$) and \reff{lipsct1},
\beq \label{corjumpestimate1}
 & & \Big|\Phi(\P_{X_r}) - \Phi(\P_{X_{r-}}) \Big|
\leq  C \E[|\Delta X_r|], \;\;\; \Big|\frac{\delta \Phi}{\delta \mu}(\P_{X_{r}}, X_{r}) - \frac{\delta \Phi}{\delta \mu} (\P_{X_{r}}, X_{r-})\Big| \leq C|\Delta X_r|.
\label{corjumpestimate2}
\enq
The claim is now clear by \reff{Xassum} from assumption {\bf (H)} and \reff{corjumpestimate2}.
\end{Remark}

\vspace{1mm}

Under such alternative conditions on the functional $\Phi$ and the semimartingale $X$, one can still derive the same form of It\^o's formula \reff{itoPP}.
For instance, one can relax the condition on $\Phi$ and instead impose a slightly stronger integrability condition on the semimartingale $X$.

\vspace{1mm}

\noindent  {\bf (H)$_{strict}$} \; There exists one decomposition $X$ $=$ $X(0)$ $+$ $V$ $+$ $L$, where $V$ is c\`adl\`ag process of finite variation with $V(0)$ $=$ $0$, and $L$ is  a local martingale with $L(0)$ $=$ $0$ such that
\beq \label{strongerX}
\E\big[\big|{\rm Var}(V)_T\big|^2\big] < \infty, \;\;\;  \E\Big[\Big(\sum_{0 < t \leq T} |\Delta X_t|\Big)^2\Big] < \infty, \;\;\; \E\big[\big|[X, X]_T^c\big|^2\big] < \infty,
\enq
where  ${\rm Var}(V)_T$, $X_{t-}$, $\Delta X_t$,  $[X, X]$ and $[X, X]^c_t$ are given in {\bf (H)}.

\begin{Theorem}[It\^o's formula (II)]\label{Theoremito_diffcond} Given a semimartingale $X$ satisfying  \bf (H)$_{strict}$
and a functional $\Phi: \Pc_2(\R^d) \to \R$ in $\Cc^{1, 1}_{relax}(\Pc_2(\R^d))$,
It\^o's formula \reff{itoPP} holds.
\end{Theorem}

In the special case when  the semimartingale $X$ is continuous, we have the following
corollary, recovering earlier results from \cite{CD2018I} and \cite{CDD2014}.

\begin{Corollary} \label{corolcont}
Suppose that the semimartingale $X$ is continuous (and hence $\P_{X_t}$ is continuous in time), and that $X$ satisfies
\beqs
\E\big[|{\rm Var}(V)_T|^2\big] < \infty,  \;\;\; \E\big[|[X, X]_T^c|^2\big] < \infty.
\enqs
Let  $\Phi$ be a function from  $\Pc_2(\R^d)$ into $\R$ satisfying
\beqs
\sup_{\mu \in \Kc} \biggl[\int_{\R^d} |\partial_x \partial_\mu \Phi(\mu, x)|^2 \mu(dx)\bigg] < + \infty,
\enqs
for any compact set $\Kc \subset \Pc_2(\R^d)$.
Then, for all $0\leq t < s\leq T$,
\beqs
\Phi(\P_{X_s }) - \Phi(\P_{X_t})
& = &
  \E \biggl[ \int_t^s \partial_\mu \Phi (\P_{X_{r}}, X_{r}).d X_r + \frac{1}{2}{\rm Tr}\Big(\partial_x \partial_\mu \Phi(\P_{X_{r}}, X_{r}) d[ X,  X]_r^c \Big)\biggl]. \nonumber
\enqs
\end{Corollary}

\vspace{1mm}

To avoid confusion, for the rest of the paper, we will always refer to the first version of It\^o's formula in Theorem \ref{TheoremIto} unless otherwise specified.

To extend It\^o's formula to the time-space-measure-dependent case, let us first define the space  $\Cc^{1, 2, (1, 1)}([0, T] \times \R^d \times \Pc_2(\R^d))$
of continuous functions $\Phi$ on  $[0, T] \times \R^d \times \Pc_2(\R^d)$ such that
\begin{itemize}
 \item  $\partial_t\Phi(t, x, \mu)$, $\partial_x\Phi(t, x, \mu)$ and $\partial_{xx}\Phi(t, x, \mu)$ exist and are jointly continuous with respect to $(t, x, \mu) \in [0, T] \times \R^d \times \Pc_2(\R^d) $;
\item both $\partial_\mu\Phi(t,x,\mu,x')$ and $\partial_{x'}\partial_\mu\Phi(t, x,\mu,x')$ exist and are jointly continuous with respect to $(t, x, \mu, x') \in [0, T] \times \R^d \times \Pc_2(\R^d)  \times \R^d$; moreover, there exists $C>0$ such that $|\partial_\mu\Phi(t, x, \mu, x')| \leq C$ and  $|\partial_{x'}\partial_\mu\Phi(t, x, \mu, x')|\leq C$.  
\end{itemize}

\begin{Corollary}[Time-space-measure-dependent It\^o's formula]\label{corollary2}
Let $X$ be a semimartingale satisfying  {\bf (H)}  (with $t$ $\mapsto$ $\P_{X_t}$ not necessarily continuous).
For any $\Phi \in \mathcal{C}^{1, 2, (1,1)}([0,T]\times \mathbb{R}^d \times \mathcal{P}_2(\mathbb{R}^d))$,
we have for $0\leq t < s\leq T$:
\beqs
 & & \Phi(s,X_s ,\P_{X_s }) =  \Phi(t, X_t, \P_{X_t})  + \int_t^s \partial_r \Phi(r,X_r ,\P_{X_{r } })dr \nonumber\\
& & +\;
   {\bar{\E}}\biggl  [ \int_t^s \partial_\mu \Phi(r,X_{r-} ,\P_{X_{r-} },{\bar X}_{r-}). d{\bar X}_{r} +  \frac{1}{2} {\rm Tr}\Big(\partial_{x} \partial_\mu  \Phi (r,X_{r-} ,\P_{X_{r-} },{\bar X}_{r-} ) d[{\bar X} ,{\bar X}]_r^c \Big)\nonumber\\
&  &+ \; \sum_{t < r \leq s} \Big\lbrace \Big(\Phi(r,X_{r-},\P_{X_{r} })  - \Phi(r,X_{r-} ,\P_{X_{r-} })\Big)  1_{\{\P_{X_r} \neq \P_{X_{r-}}\}}
  -   \partial_\mu \Phi (r,X_{r-} ,\P_{X_{r-} },{\bar X}_{r-} ). \Delta {\bar X}_{r}\Big\rbrace\biggl] \nonumber\\
& & + \; \bar\E\Big[ \sum_{t < r \leq s} \biggl (\frac{\delta \Phi}{\delta \mu}(r,X_{r-} ,\P_{X_{r}}, \bar X_{r}) -  \frac{\delta \Phi}{\delta \mu}(r,X_{r-} ,\P_{X_{r-}}, \bar X_{r-}) \biggl) 1_{\{\P_{X_r} = \P_{X_{r-}}\}} \Big] \nonumber\\
& & +  \; \sum_{i =1}^d  \int_t^s \partial_{x_i} \Phi (r,X_{r-} ,\P_{X_{r } }) dX_{r}^i
  +   \frac{1}{2}\sum_{i,j=1}^d  \int_t^s  \partial_{x_i}\partial_{x_j} \Phi  (r,X_{r-} ,\P_{X_{r } }) d[X^i ,X^j ]_r^c \nonumber\\
& &+\; \sum_{t < r \leq s} \biggl \lbrace \Phi(r,X_{r } ,\P_{X_r}) - \Phi(r,X_{r-} ,\P_{X_r})
 -\sum_{i =1}^d \partial_{x_i}\Phi (r,X_{r-} ,\P_{X_r}) \Delta X_{r}^i
 \biggl \rbrace.
  \label{itoTXP}
 \enqs
 Here $\bar X$ denotes an independent and identical copy of $X$ on a copy of the probability space $(\Omega, \Fc, \P)$.
\end{Corollary}
{\bf Proof.} To see how it follows from Theorem \ref{TheoremIto},
fix an $m$ $\in$ $\N$. Let $\pi^m_{t, s}$ $=$ $(t=t_0^m < t_1^m< \cdots < t_{m + 1}^m = s)$ be any partition of $[t,s]$.
\begin{align*}
& \Phi(s, X_s ,\P_{X_s}) - \Phi(t, X_t, \P_{X_t}) \\
&= \sum_{j=0}^{m} \biggl \lbrace \Phi(t_{j+1}^m, X_{t_{j+1}^m} ,\P_{X_{t_{j+1}^m} }) - \Phi(t_j^m, X_{t_j^m} ,\P_{X_{t_j^m} }) \biggl \rbrace
\\ &= \sum_{j=0}^{m} \biggl \lbrace \Phi(t_{j+1}^m, X_{t_{j+1}^m} ,\P_{X_{t_{j+1}^m} })-\Phi(t_j^m, X_{t_j^m},\P_{X_{t_{j+1}^m} }) + \Phi(t_j^m, X_{t_j^m} ,\P_{X_{t_{j+1}^m} }) - \Phi(t_j^m, X_{t_j^m} ,\P_{X_{t_j^m} }) \biggl \rbrace
\\  &= \sum_{j=0}^{m} \biggl \lbrace \Phi(t_{j+1}^m, X_{t_{j+1}^m} ,\P_{X_{t_{j+1}^m}})-\Phi(t_j^m, X_{t_j^m} ,\P_{X_{t_{j+1}^m}}) \biggl \rbrace +  \sum_{j=0}^{m} \biggl \lbrace \Phi(t_j^m, X_{t_j^m},\P_{X_{t_{j+1}^m}}) - \Phi(t_j^m, X_{t_j^m},\P_{X_{t_j^m}}) \biggl \rbrace.
\end{align*}
Let $m \to \infty$, then, similar as proving the classical It\^o's formula for semimartingales, we see that the first sum converges in probability to
\begin{align*}
&  \int_t^s \partial_r \Phi(r,X_r ,\P_{X_{r } })dr + \sum_{i=1}^d\int_t^s \partial_{x_i} \Phi (r, X_{r-} ,\P_{X_{r} }) dX_{r}^i
  +   \frac{1}{2} \sum_{i, j=1}^d \int_t^s  \partial_{x_ix_j}^2 \Phi  (r, X_{r-} ,\P_{X_{r } }) d[X^i,X^j]_r^c
 \\&+ \sum_{t < r \leq s}\biggl \lbrace \Phi(r, X_{r } ,\P_{X_{r } }) - \Phi(r, X_{r-},\P_{X_{r } })
  -\sum_{i=1}^d \partial_{x_i} \Phi (r, X_{r-} ,\P_{X_{r } }) \Delta X_{r}^i \biggl \rbrace.
\end{align*}
According to Theorem \ref{TheoremIto}, the second sum converges to
\begin{align*}
& { \bar {\E}}\biggl  [ \int_t^s \partial_\mu \Phi(r, X_{r-} ,\P_{X_{r-} }, \bar {X}_{r-} ). d{\bar X}_{r}
  +  \int_t^s {\rm Tr} \Big(\partial_x\partial_\mu  \Phi (r, X_{r-} ,\P_{X_{r-} }, \bar {X}_{r-} ) d[{\bar X} ,{\bar X} ]_r^c\Big) \biggl ]
\\ &+ \sum_{t < r \leq s} \Big ( \Phi(r, X_{r-} ,\P_{X_{r} })  - \Phi(r, X_{r-} ,\P_{X_{r-}})\Big ) 1_{\{\P_{X_r} \neq \P_{X_{r-}}\}}
 -   {\bar {\E}} \Big [\sum_{t < r \leq s} \partial_\mu \Phi (r, X_{r-},\P_{X_{r-}}, \bar {X}_{r-}) \cdot \Delta { \bar X}_{r}\Big ] \\
 &  +  \bar\E\Big[\sum_{t < r \leq s} \Big (\frac{\delta \Phi}{\delta \mu}(r,X_{r-} ,\P_{X_{r-}}, \bar X_{r}) -  \frac{\delta \Phi}{\delta \mu}(r,X_{r-} ,\P_{X_{r-}}, \bar X_{r-}) \Big) 1_{\{\P_{X_r} = \P_{X_{r-}}\}}  \Big]\nonumber.
\end{align*}
\ep

\vspace{1mm}

It\^o's formula in Theorem \ref{TheoremIto} takes some special forms when applied  to a  class of semimartingales driven by the McKean--Vlasov SDEs or jump diffusion processes which are useful for the subsequent analysis of McKean--Vlasov controls in Sections  \ref{secMKVjump} and \ref{secMKVsingular}.

To see this, first let $\boldmath{b}=(b_s)_{s \geq 0}$ and $\boldmath{\sigma}=(\sigma_s)_{s \geq 0}$ be $\F$-adapted processes valued in $\R^d$ and $\R^{d \times d}$, respectively, and $\beta = (\beta_s(\theta))_{s \geq 0}$ be $\F$-predictable processes valued in $\R^d$ for $\theta \in \R^q$. Consider the following jump-diffusion process $(X_s)_{s \geq 0}$
\beq \label{equMKVSDEjump}
d X_s = b_s ds + \sigma_s dW_s + \int_{\theta \in \R^q} \beta_s(\theta) N(ds, d\theta),\;\;\; X_{t-} = \xi \in L^2(\Fc_t; \R^d),
\enq
where $(W_s)_{s \geq 0}$ is Brownian motion, and $N(ds, d\theta)$ is a Poisson random measure with a finite intensity measure $\nu$, with its compensated Poisson random measure  $\tilde N(ds, d\theta): = N(ds, d\theta) - \nu(d\theta) ds$.  It is easy to check that such $(X_s)_{s \geq t}$ satisfies assumption ({\bf H}) if
$$ \E\Big[\int_0^T \Big(|b_s|^2 + |\sigma_s|^2 + \int_{\R^q}|\beta_s(\theta)| \nu(d\theta)\Big) ds\Big] < \infty.$$
Note that such process $X$ has jumps induced by the Poisson random measure, whereas $\P_{X}$ is continuous in time. (See, for instance \cite{G2017}).

\begin{Corollary} [It\^o's formula for jump process] \label{CorItojump}
Given a semimartingale $X$ defined by \reff{equMKVSDEjump}, then for any $\Phi \in \mathcal{C}^{1,1}(\mathcal{P}_2(\mathbb{R}^d))$, we have
\beqs
\Phi(\P_{X_s}) &=& \Phi(\P_{\xi}) + \int_t^s  \E\biggl[ \partial_\mu \Phi(\P_{X_r}, X_r).  b_r + \frac{1}{2} {\rm Tr}\Big(\sigma_r \sigma_r\trans\partial_x\partial_\mu \Phi(\P_{X_r},  X_r)\Big)\\
& & \quad \quad \quad \quad   +\; \int_{\R^q} \Big( \frac{\delta  \Phi}{\delta \mu}(\P_{X_{r}}, X_{r} + \beta_r(\theta)) - \frac{\delta  \Phi}{\delta \mu}(\P_{X_{r}},  X_{r}) \Big) \nu(d\theta)\biggl] dr,
\enqs
for all $0\leq  t < s \leq T$.
\end{Corollary}
{\bf Proof.} \  \ Apply Theorem \ref{TheoremIto} to $\Phi(\P_{X_r})$ with $X_r$ given in \reff{equMKVSDEjump} between $t$ and $s$,
\beq \label{proofItojump}
\Phi(\P_{X_{s}}) - \Phi(\P_{X_{t-}}) &=& \E\biggl[\int_t^s \partial_\mu \Phi (\P_{X_{r}}, X_{r}). b_r
+ \frac{1}{2}{\rm Tr}\Big(\sigma_r \sigma _r \trans \partial_x \partial_\mu \Phi(\P_{X_{r}}, X_{r})\Big) dr \nonumber\\
&  & \; + \;\;\; \sum_{t < r \leq s} \Big\lbrace \Big(\Phi(\P_{X_r}) - \Phi(\P_{X_{r-}})\Big) 1_{\{\P_{X_r} \neq \P_{X_{r-}}\}} \nonumber\\
& &\;+  \Big(\frac{\delta \Phi}{\delta \mu}(\P_{X_{r}}, X_{r}) - \frac{\delta \Phi}{\delta \mu} (\P_{X_{r}}, X_{r-})\Big) 1_{\{\P_{X_r} = \P_{X_{r-}}\}}  \Big\rbrace\biggl].
\enq
Now it suffices to compute the jump term in \reff{proofItojump}.
As the distribution $\P_{X_r}$ is continuous in time, we have
\beqs
\sum_{t < r \leq s} \Big\lbrace \Big(\Phi(\P_{X_r}) - \Phi(\P_{X_{r-}})\Big) 1_{\{\P_{X_r} \neq \P_{X_{r-}}\}} \Big\rbrace =0.
\enqs
Note that at the time $r$, $X_{r-}$ has a jump $\beta_r(\theta)$ with Poisson random measure $N(dr, d\theta)$, 
Thus,
\beqs
& & \E\biggl[\sum_{t < r \leq s} \Big\lbrace \Big(\Phi(\P_{X_r}) - \Phi(\P_{X_{r-}})\Big) 1_{\{\P_{X_r} \neq \P_{X_{r-}}\}} \\
& & \;+\; \Big(\frac{\delta \Phi}{\delta \mu}(\P_{X_{r}}, X_{r}) - \frac{\delta \Phi}{\delta \mu} (\P_{X_{r}}, X_{r-})\Big) 1_{\{\P_{X_r} = \P_{X_{r-}}\}} \Big\rbrace \biggl]\\
& & = \E\biggl[\int_t^{s}\int_{\theta \in \R^q} \Big( \frac{\delta  \Phi}{\delta \mu}(\P_{X_{r}}, X_{r} + \beta_r(\theta)) - \frac{\delta  \Phi}{\delta \mu}(\P_{X_{r}}, X_{r}) \Big) \nu(d\theta) dr\biggl].
\enqs
\ep

\paragraph{Relation to Fokker--Planck equation} One can check that this form of It\^o's formula for \reff{equMKVSDEjump} is consistent with its Fokker--Planck equation. To see this, suppose that $\P_{X_t}$ has a PDF $p_t$ such that $b_t = b(X_t, p_t)$, $\sigma_t = \sigma(X_t, p_t)$, $\beta_t(\theta) = \theta$, and $N$ is a compound Poisson Process with intensity $\lambda \in \R^d$ and jump PDF $\gamma$, then under mild conditions, the time evolution of $p_t$ can be prescribed by an integral differential equation
\beq \label{FPPDF}
\partial_t p_t &=& -\nabla \cdot (p_t b(x, p_t)) + \frac{1}{2} \nabla \cdot \nabla \cdot (\sigma(x, p_t)\sigma(x, p_t)\trans p_t) \nonumber\\
& & \hspace{2cm} + \; \lambda \Big(\int_{\R^d} p_t(\theta) \gamma (x-\theta)d\theta - p_t(x)\Big),
\enq
where $\nabla \cdot$ denotes the divergence. One can further obtain the time evolution of $\phi(p_t)$ for a smooth function $\phi: L^2(\R^d) \to \R$ with respect to the G\^ateaux derivative (\cite{BJS2017}),
\beqs
d\phi(p_t) &=& \int_{\R^d} \frac{\partial \phi}{\partial p}(p_t, x) \partial_t p_t dx dt \\
& =&\E\Big[\partial_x \frac{\partial \phi}{\partial p}(p_t, x). b(X_t, p_t) + \frac{1}{2} {\rm Tr}\Big(\sigma\sigma\trans(X_t, p_t) \partial_{xx} \frac{\partial \phi}{\partial p}(p_t, X_t)\Big)\\
& & \; + \; \lambda\int_{\R^d} \Big(\frac{\partial \phi}{\partial p}(p_t, X_t + \theta) - \frac{\partial \phi}{\partial p}(X_t)\Big) \gamma(\theta) d\theta \Big] dt,
\enqs
where the second inequality is by \reff{FPPDF} and from integration by parts. Now the consistency follows from  the connection \eqref{reldermu}  between the G\^ateaux derivative and the linear derivative in Definition \ref{linear}.

 Next, let us consider processes that arise from  control problems of  singular type. Let $(\eta_t)$ $=$ $(\eta_t^i)_{1 \leq i \leq d}$ be finite variation processes valued in $\R^d$, and denote $\lambda$ $=$ ${\rm diag}(\lambda_1, \ldots, \lambda_d)$ with
$\lambda_i$ nonnegative constant in $\R$. Take the following c\`adl\`ag process $X$ $=$ $(X_s^1, \ldots, X_s^d)_{s \geq t}$,
\beq \label{XSDE}
 dX_s = b_s ds + \sigma_s dW_s +   \lambda d \eta_s, \; s \geq t, \; X_{t-} = \xi
 \in L^2(\Fc_t;\R^d),
 \enq
 where $(W_s)_{s \geq 0}$  is a $d$-dimensional Brownian motion. Assuming that  $\E[{\rm Var}(\eta)_T] < \infty$ and $\E[\int_0^T (|b_s|^2 + |\sigma_s|^2) ds] < \infty $,
 one  can easily check that $X$ satisfies  assumption {\bf (H)}, and  we then have,  

\begin{Corollary} \label{CorItosingular}
Given a semimartingale $X$ in \reff{XSDE} (with $t \mapsto \P_{X_t}$ not necessarily continuous). For any $\Phi \in \mathcal{C}^{1,1}(  \mathcal{P}_2(\mathbb{R}^d))$,
\beqs\label{itofor}
& &\Phi( \P_{X_s}) - \Phi(\P_\xi) \nonumber\\
& = &   {\mathbb{E}} \biggl [ \int_t^s\Big \lbrace\partial_\mu \Phi (\P_{X_{r}},{X}_{r}).b_r + \frac{1}{2} {\rm Tr}
\Big( \partial_{x} \partial_\mu  \Phi (\P_{X_{r}}, {X}_{r}) \sigma_r \sigma_r\trans\Big) \Big\rbrace dr \nonumber\\
& &  +  \; \sum_{i=1}^d \limits  \int_t^s \lambda_i \partial_\mu \Phi_i (\P_{X_{r-}}, {X}_{r-})  d{\eta}_r^i \nonumber\\
& & + \;\sum_{t \leq r \leq s} \Big \lbrace \Big(\Phi( \P_{X_r^{}})  - \Phi( \P_{X_{r-}^{}})\Big) 1_{\{\P_{X_r} \neq \P_{X_{r-}}\}} -\sum_{i=1}^d \lambda_i \partial_\mu \Phi_i (\P_{X_{r-}^{}},{X}_{r-}^{})  \Delta{\eta}_r^i \Big\rbrace \biggl] \nonumber\\
&  & + \;  \E\biggl[\sum_{t \leq r \leq s}\Big(\frac{\delta \Phi}{\delta \mu}(\P_{X_{r}}, X_{r}) - \frac{\delta \Phi}{\delta \mu} (\P_{X_{r}}, X_{r-})\Big) 1_{\{\P_{X_r} = \P_{X_{r-}}\}} \biggl], \quad \quad 0 \leq t \leq s \leq T.
\enqs
\end{Corollary}

\subsection{Proof of Theorem \ref{TheoremIto}}

We will establish It\^o's formula  for a class of cylindrical functions called $\Gc(\Mc)$. We will then prove the general case by applying  the general form of Stone--Weierstrass theorem on compact sets in the joint Wasserstein space and Euclidean space and by appropriate localization argument.  For ease of exposition, we will prove for the case of $d$ $=$ $1$ without loss of generality. We will discuss (see Remark \ref{remd1}) how to adapt the arguments to the general case of $d$ $>$ $1$.

Throughout the proof of Theorem \ref{TheoremIto}, we adopt a generic constant $C$ for ease of exposition, unless otherwise specified.

\subsubsection{Proof of Theorem \ref{TheoremIto} for  $\Gc(\Kc_K)$}

Let us first  define the general space of cylindrical functions $\Gc(\Mc)$ for  $\Mc \subset \Pc_2(\R)$.

\begin{Definition} \label{defGM} Given $\Mc$ $\subset$ $\Pc_2(\R)$, define
\begin{align} \label{GcM}
\Gc(\Mc) :=\{& \Mc \ni \mu  \mapsto  \Phi ( \mu) = f\big (\left<g_1, \mu\right>,\cdots,\left<g_n, \mu\right>\big) \text{ for some } n \in \N, \nonumber\\ & \text{and polynomial } f :\mathbb{R}^n  \rightarrow \mathbb{R} \text{ and } g_1,\cdots, g_n  :\mathbb{R} \rightarrow \mathbb{R} \text{ polynomials }\}.
\end{align}
\end{Definition}
Observe that $\Gc(\Mc)$ is an algebra as it is closed under pointwise addition, multiplication, and scalar multiplication. Moreover, $\Gc(\Mc)$ can be rewritten in the following form
\begin{align*}
\Gc(\Mc) = \mathrm{span}\big\lbrace{\left<g, \mu^n\right>: n \geq 1, \mu \in \Mc, g: \R^n \to \R \text{ is monomial}\big\rbrace},
\end{align*}
where $\left<g, \mu^n\right>: = $ $\int_{\R^k} g(x_1, \ldots, x_n)\mu(dx_1)\ldots \mu(dx_n)$ with $n$ $\geq$ $1$ is a basis of $\Gc(\Mc)$, as each $\Phi \in \Gc(\Mc)$ is a linear combination of monomials $\left<g, \mu^k\right>$. (See \cite{CLSF2019} for more discussions of this function space).

\medskip

Now,  consider the cylindrical function $\Phi(\mu) = f(\left<g_1, \mu\right>, \ldots, \left<g_n, \mu \right>)$, whose Lions derivative is given by
\begin{align*}
\partial_\mu \Phi(\mu, x) &= \;
\sum_{k=1}^n\partial_{y_k} f(\left<g_1, \mu\right>, \ldots, \left<g_n, \mu\right>) \partial_{x} g_k(x).
\end{align*}
Note, however, the Lions derivative and the mixed second order derivative of cylindrical functions are not necessarily in $\Gc(\Mc)$. Instead they belong to a bigger algebraic space defined below.
\begin{Definition}\label{defHMU}
Given $\Mc$ $\subset$ $\Pc_2(\R)$ and $U$ $\subset$ $\R$. Let
\begin{align*}
\Hc( \Mc \times U) :=\{&\Mc \times U \ni (\mu, x) \mapsto \Phi(\mu, x) =  \sum_{k=1}^n f_k\big (\left<g_k, \mu^k\right>\big) h_k(x),  \text{ for some } n \in \N,\\ & f_k, h_k:\mathbb{R} \rightarrow \mathbb{R} \text{ and } g_k: \R^k \to \R \text{ are monomials }\}.
\end{align*}
\end{Definition}
One can see that $\Hc(\Mc \times U)$ is also an algebra and can be rewritten in the following form
\begin{align*}
\Hc(\Mc \times U) = \mathrm{span} \biggl\lbrace{\left<g, \mu^k\right> x^l: k \geq 1, l \geq 0, \mu \in \Mc, x \in U, g: \R^k \to \R \text{ is monomial} \biggl \rbrace}.
\end{align*}
Moreover, $\Gc(\Mc)$ can be viewed as a subalgebra of $\Hc(\Mc \times U)$.

\medskip

We will next consider a particular choice of  $\Gc(\Mc)$ with $\Mc:=\Kc_K$, and
\beqs
\Kc_K: = \Big\{\mu  \in \Pc_2(\R)| \mathrm{supp}(\mu) \subset [-K, K]\Big\}.
\enqs
Clearly $\Kc_K$ is not empty since the Dirac measure $\delta_0$ is in $\Kc_K$.
Moreover,  by Lemma 5.7 and Proposition 5.3 in \cite{C2013}, $\Kc_K$ is compact under the $W_p$ distance for any $p$ $\geq$ $1$. Now we can establish It\^o's formula for any $\Phi \in \Gc(\Kc_K)$.

\begin{Lemma}\label{M2C}
Given a semimartingale $X$ satisfying assumption {\bf (H)} and $|X_t| \leq K$ $\P$-a.s., for some $K > 0$ and any $t \in [0, T]$.
Then It\^o's lemma in the form of equation \reff{itoPP} holds for $\Phi(\P_{X_t})$, with $\Phi$ $\in$ $\Gc(\Kc_K)$.
\end{Lemma}
{\bf Proof.}\;
Given $\Phi \in \Gc(\Kc_K)$, $\Phi(\mu) = f\big (\left<g_1, \mu\right>,\cdots,\left<g_n, \mu\right>\big)$ by \reff{GcM}, where $f$ and $g_k$, $1 \leq k \leq n$  are polynomials. Therefore,
\beqs
\Phi(\P_{X_r }) = f(\mathbb{E}[ g_1(X_r)], \ldots, \E[g_n(X_r)]),
\enqs
where the process $r \mapsto X_r$ is a semimartingale and $|X_r|$ $\leq$ $K$ $\P$-a.s.
Clearly, if we define $Y_r^k =  \mathbb{E}[ g_k(X_r)]$, $k =1, \ldots, n$ for any $r \in [0, T]$, then $ \Phi(\P_{X_r }) = f(Y_r^1, \ldots, Y_r^n)$. Since $X_r$ is a semimartingale and $g_k$ is polynomial, by the classical It\^o's formula for semimartingale, $g_k(X_r)$ and (hence) $Y_r^k = \mathbb{E}[g_k(X_r)]$ are  bounded semimartingales.
Now, setting $g$ $=$ $(g_1,\ldots,g_n)$ and $Y_r$ $=$ $(Y^1_r,\ldots,Y^n_r)$, $0 \leq r \leq T$, and
applying It\^o's formula to $f(Y_r)$ between $t$ and $s$, $0$ $\leq$ $t$ $<$ $s$ $\leq$ $T$, we see
\begin{align} \label{proofItoPhi}
 &\Phi(\P_{X_s }) - \Phi(\P_{X_t}) =  f(Y_s) - f( Y_t)\nonumber \\
 &=  \underbrace{\sum_{k=1}^n \int_t^s \partial_{y_k} f( Y_{r-}) dY_r^k}_{I}
 + \underbrace{\frac{1}{2} \sum_{k=1}^n \sum_{j=1}^n\int_t^s \partial_{y_ky_j}  f ( Y_{r-}) d[Y^j,Y^k]_r^c}_{II} \\
 &+ \underbrace{\sum_{t < r \leq s} \biggl\lbrace  f( Y_{r}) -  f(Y_{r-}) - \sum_{k=1}^n \partial_{y} f( Y_{r-}^k) \Delta Y_r^k \biggl\rbrace}_{III}. \nonumber
\end{align}

Let us now compute the terms $I$, $II$, and $III$ separately.

\vspace{1mm}

Applying the classical It\^o's formula to $g_k(X_r)$, $1 \leq k \leq n$, we have
\beq \label{equ:Itogk}
g_k(X_s) - g_k(X_t) &=& \int_t^s \partial_x g_k(X_{r-}) dX_r + \frac{1}{2} \partial_{xx} g_k(X_{r-}) d[X, X]_r^c\\
& & \;\;\; +\;  \sum_{t < r \leq s} \Big\lbrace g_k(X_r) - g_k(X_{r-}) - \partial_x g_k(X_{r-}) \Delta X_r\Big\rbrace. \nonumber
\enq
Since $|X_r|$ $\leq$ $K$, $\P$-a.s., there exists constant $C$ $>$ $0$ such that $\P$-a.s.
\beqs
\Max\Big(|g_k(X_r)|, |\partial_x g_k(X_{r-})|, |\partial_{xx} g_k(X_{r-})|\Big) \leq C.
\enqs
Taking expectation on both sides of \reff{equ:Itogk}, we see
\beq \label{YIto}
Y_s^k - Y_t^k &=& \E\biggl[\int_t^s \partial_x g_k(X_{r-}) d X_r + \frac{1}{2} \partial_{xx} g_k(X_{r-}) d[X, X]_r^c \\
& & \; + \; \sum_{t < r \leq s} \Big\lbrace g_k(X_r) - g_k(X_{r-}) - \partial_x g_k(X_{r-}) \Delta X_r\Big\rbrace \biggl]. \nonumber
\enq
By the definition of the integral,
\beqs \label{stoY}
I: = \sum_{k=1}^n\int_t^s \partial_{y_k} f( Y_{r-}) dY_r^k = \lim_{m \to \infty}\sum_{\pi_{t, s}} \sum_{k=1}^n \partial_{y_k} f(Y_{t_j^m}) (Y_{t_{j+1}^m}^k - Y_{t_j^m}^k),
\enqs
holds for any arbitrary partition $\pi^m_{t, s}$ of $[t, s]$.
Therefore, in order to calculate the term $I$, we may take a partition $\bar\pi^{m}_{t, s}$ $=$ $(t=t_0^m < t_1^m< \cdots < t_{m + 1}^m = s), m \in \N$ such that $\max_j|t_j^m - t_{j+1}^m| \to 0$ and  $\max_{j} \sup_{r \in (t_j^m, t_{j+1}^m]}|Y_{r-} - Y_{t_j^m}| \to 0$ as $m \to \infty$. The existence of $\bar\pi^{m}_{t, s}$ is ensured by the fact that $(Y_r)_{t \leq r \leq  s}$ is a c\`adl\`ag function of finite variation. Therefore, $I$ is well-defined and can be approximated by the partition
$\bar\pi^{m}_{t, s}$:
\beqs
I&=&  \lim_{m \to \infty}\sum_{\bar\pi_{t, s}} \sum_{k=1}^n \partial_{y_k} f(Y_{t_j^m}) (Y_{t_{j+1}^m}^k - Y_{t_j^m}^k) \\
&=& \lim_{m \to \infty}\sum_{j=0}^{m} \E\biggl[ \int_{t_j^m}^{t_{j+1}^{m}}\sum_{k=1}^n \partial_{y_k} f(Y_{t_j^m})\partial_x g_k(X_{r-}) dX_r\\
& & \hspace{1cm}+\; \frac{1}{2} \int_{t_j^m}^{t_{j+1}^{m}}\sum_{k=1}^n \partial_{y_k} f(Y_{t_j^m})\partial_{xx} g_k(X_{r-}) d[X, X]_r^c \nonumber\\
& & \hspace{1cm}+\; \sum_{t_j^m < r \leq t_{j+1}^m} \sum_{k=1}^n \partial_{y_k} f(Y_{t_j^m})\Big\lbrace g_k(X_r) -g_k(X_{r-}) - \partial_x g_k(X_{r-}) \Delta X_r\Big\rbrace\biggl] \nonumber\\
&=& \E\biggl[ \int_t^s (\partial_\mu \Phi )(\P_{X_{r-} }, {X}_{r-} ) d{X}_{r} + \frac{1}{2}\partial_{x}\partial_\mu \Phi(\P_{X_{r-} },{X}_{r-})d[{X} ,{X}]_r^c \\
& & + \;\sum_{t < r \leq s} \biggl\lbrace \frac{\delta \Phi}{\delta \mu}(\P_{X_{r-}}, X_r) - \frac{\delta \Phi}{\delta \mu}(\P_{X_{r-}}, X_{r-}) - \partial_\mu \Phi(\P_{X_{r-}}, X_{r-}) \Delta X_{r} \biggl\rbrace\biggl],
\enqs
where in the second equality, \reff{YIto} is applied to $Y^k_t$ between $t_j^m$ and $t_{j + 1}^m$, and in the last equality, $Y_{t_j^m}$ in the second inequality is replaced by $Y_{r-}$, $r \in (t_j^m, t_{j+1}^m]$ as $\max_{j} \sup_{r \in (t_j^m, t_{j+1}^m]}|Y_{r-} - Y_{t_j^m}| \to 0$ as $m \to \infty$.
\medskip

Now, since $Y^k$ is deterministic, $[Y^k, Y^j]_r^c =0$, and
\beqs
II: = \frac{1}{2} \sum_{k=1}^n\sum_{j=1}^n \int_t^s \partial_{y_ky_j}  f (Y_{r-}) d[Y^k,Y^j]_r^c =0.
\enqs

\medskip
Moreover,  notice that
the third term of \reff{proofItoPhi} is
\beqs
III & & = \sum_{t < r \leq s} \biggl\lbrace \Phi(\P_{X_r}) - \Phi(\P_{X_{r-}}) - \sum_{k=1}^n  \E\Big[\partial_{y_k}
f(Y_{r-})\big(g_k(X_r) - g_k( X_{r-})\big)\Big] \biggl\rbrace\\
& & = \sum_{t < r \leq s} \biggl\lbrace \Phi(\P_{X_r}) - \Phi(\P_{X_{r-}}) - \E\Big[\frac{\delta \Phi}{\delta \mu}(\P_{X_{r-}}, X_r) - \frac{\delta \Phi}{\delta \mu}(\P_{X_{r-}}, X_{r-})\Big] \biggl\rbrace.
\enqs
Summing up these three terms $I$, $II$, and $III$, and substituting them into \reff{proofItoPhi}, we obtain
\begin{align*}
 & \Phi(\P_{X_s }) - \Phi(\P_{X_t}) \\
 = & \;\;
   {\mathbb{E}} \biggl[ \int_t^s  \partial_\mu \Phi (\P_{X_{r-} },{X}_{r-} ) d{X}_{r} + \frac{1}{2} \partial_{x}\partial_\mu \Phi (\P_{X_{r-} },{X}_{r-})  d[{X}  ,{X} ]_r^c \biggl]\\
 & \; + \sum_{t < r \leq s} \Big(\Phi(\P_{X_{r}} ) - \Phi(\P_{X_{r-}})\Big) 1_{\{\P_{X_r} \neq \P_{X_{r-}}\}}-   {\mathbb{E}} \Big[ \sum_{t < r \leq s} \partial_\mu \Phi   (\P_{X_{r-} },{X}_{r-} ) \Delta {X}_{r} \Big]\\
 & \; + \E\biggl[\sum_{t < r \leq s} \biggl\lbrace \frac{\delta \Phi}{\delta \mu}(\P_{X_{r-}}, X_r) - \frac{\delta \Phi}{\delta \mu}(\P_{X_{r-}}, X_{r-})\biggl\rbrace\biggl]\\
& \; - \sum_{t < r \leq s}\E\Big[\Big(\frac{\delta \Phi}{\delta \mu}(\P_{X_{r-}}, X_{r}) - \frac{\delta \Phi}{\delta \mu} (\P_{X_{r-}}, X_{r-})\Big)\Big].
\end{align*}
As  $(\P_{X_r})_{t < r \leq s }$ is a c\`adl\`ag function, let us divide $[t, s]$ into countable set $\{r \in [t, s]: \P_{X_r} \neq \P_{X_{r-}}\}$ and uncountable set  $\{r \in [t, s]: \P_{X_r} = \P_{X_{r-}}\}$, we have
\beqs
& & \E\biggl[\sum_{t < r \leq s} \biggl(\frac{\delta \Phi}{\delta \mu}(\P_{X_{r-}}, X_r) - \frac{\delta \Phi}{\delta \mu}(\P_{X_{r-}}, X_{r-})\biggl)\biggl]\\
& = &   \sum_{t < r \leq s} \E\biggl[\biggl(\frac{\delta \Phi}{\delta \mu}(\P_{X_{r-}}, X_r) - \frac{\delta \Phi}{\delta \mu}(\P_{X_{r-}}, X_{r-})\biggl)1_{\{\P_{X_r} \neq  \P_{X_{r-}}\}}\biggl] \\
& &\;\; + \;  \E\biggl[\sum_{t < r \leq s} \biggl(\frac{\delta \Phi}{\delta \mu}(\P_{X_{r-}}, X_r) - \frac{\delta \Phi}{\delta \mu}(\P_{X_{r-}}, X_{r-}) \biggl)1_{\{\P_{X_r} = \P_{X_{r-}}\}}\biggl],
\enqs
where we switch the order of the expectation and the summation on countable set $\{r \in [t, s]: \P_{X_r} \neq \P_{X_{r-}}\}$. Meanwhile,
\beqs
& & \sum_{t < r \leq s}\E\Big[\Big(\frac{\delta \Phi}{\delta \mu}(\P_{X_{r-}}, X_{r}) - \frac{\delta \Phi}{\delta \mu} (\P_{X_{r-}}, X_{r-})\Big)\Big]\\
& = &   \sum_{t < r \leq s}\E\Big[\Big(\frac{\delta \Phi}{\delta \mu}(\P_{X_{r-}}, X_{r}) - \frac{\delta \Phi}{\delta \mu} (\P_{X_{r-}}, X_{r-})\Big)\Big] 1_{\{\P_{X_r} \neq \P_{X_{r-}}\}}.
\enqs
We conclude that
\beqs
& & \E\biggl[\sum_{t < r \leq s} \biggl\lbrace \frac{\delta \Phi}{\delta \mu}(\P_{X_{r-}}, X_r) - \frac{\delta \Phi}{\delta \mu}(\P_{X_{r-}}, X_{r-})\biggl\rbrace\biggl]
- \sum_{t < r \leq s}\E\Big[\Big(\frac{\delta \Phi}{\delta \mu}(\P_{X_{r-}}, X_{r}) - \frac{\delta \Phi}{\delta \mu} (\P_{X_{r-}}, X_{r-})\Big)\Big]\\
& = &  \E\Big[\sum_{t < r \leq s} \Big\lbrace \frac{\delta \Phi}{\delta \mu}(\P_{X_{r-}}, X_r) - \frac{\delta \Phi}{\delta \mu}(\P_{X_{r-}}, X_{r-}) \Big\rbrace 1_{\{\P_{X_r} = \P_{X_{r-}}\}} \Big].
\enqs
Note that $|\partial_\mu\Phi(\mu, x)|$ is bounded on the compact set $\Kc_K \times [-K, K]$, and by the same argument  as in Remark \ref{Itoremainder}, one can see that the terms
$ \sum_{t < r \leq s} \Big(\Phi(\P_{X_{r}} ) - \Phi(\P_{X_{r-}} )\Big)1_{\{\P_{X_r} \neq \P_{X_{r-}}\}}$,   ${\mathbb{E}} \Big[\sum_{t < r \leq s}\partial_\mu \Phi   (\P_{X_{r-} },{X}_{r-} ) \Delta {X}_{r}\Big]$, and $\E\Big[\sum_{t < r \leq s} \Big\lbrace \frac{\delta \Phi}{\delta \mu}(\P_{X_{r-}}, X_r) - \frac{\delta \Phi}{\delta \mu}(\P_{X_{r-}}, X_{r-}) \Big\rbrace 1_{\{\P_{X_r} = \P_{X_{r-}}\}} \Big]$ are all finite.
\ep

\subsubsection{Proof of Theorem \ref{TheoremIto} for
$\Cc^{1, 1}(\Kc_K)$}

We next  establish the  It\^o's lemma for $\Cc^{1, 1}(\Kc_K)$. To this end, we first need:
\begin{Lemma}[Stone-Weierstrass] \label{LemmaSW}
\label{SWlemma}
Take  a compact Hausdorff space $H$, and let $\mathcal{C}(H)$ be the algebra of real-valued continuous functions on $H$, with the topology of uniform convergence. Let  $\mathcal{A}$ be a subalgebra of  $\mathcal{C}(H)$.
If   $\mathcal{A}$ separates points on $H$ and vanishes at no point on $H$, then $\mathcal{A}$ is dense in $\mathcal{C}(H)$.
\end{Lemma}

Then we will establish  the following Lemma showing that  cylindrical functions restricted on $\Kc_K$ is dense in $\Cc^{1, 1}(\Kc_K)$ with appropriate choices of norms.
\begin{Lemma} \label{Pdensity}$\Gc(\Kc_K)$ is dense in $\Cc^{1, 1}(\Kc_K)$, the collection of all $\Cc^{1, 1}$ functions on $\Kc_K$ with the supremum norm of derivatives of all orders:
\beqs
\|\Phi\|_{\Kc_K} &: =& \sup_{(\mu, x) \in \Kc_K \times [-K, K]}\Big(|\Phi(\mu)|  + |\partial_\mu \Phi(\mu, x)| + |\partial_x\partial_\mu \Phi(\mu, x)|\Big).
\enqs
\end{Lemma}
{\bf Proof.} The proof consists of  two steps, and we will adopt the linear derivative.

\noindent {\it Step 1.}
 If $\Phi$ $\in$ $\Cc^{1, 1}(\Kc_K)$, then $\partial_{xx} \frac{\delta \Phi}{\delta \mu}(\mu, x)$ $\in$ $\Cc(\Kc_K \times [-K, K])$. Since we are concerned about the support of $\mu$, we restrict $\partial_{xx} \frac{\delta \Phi}{\delta \mu}$ to $\Kc_K$ $\times$ $[-K, K]$. We see that
\begin{itemize}
\item {\it $\Hc(\Kc_K \times [-K, K])$ $=:$ $\Hc_K$ separates points on $\Kc_K$ $\times$ $[-K, K]$}. To see this, take $(\mu, x) \neq (\mu', x') \in \mathcal{K}_K \times [-K, K]$, we have either $\mu$ $\neq$ $\mu'$ or $x$ $\neq$ $x'$. If $\mu$ $\neq$ $\mu'$, from Theorem 30.1 in \cite{B2013}, there exists $j_0$ $\in$ $\N$ such that $\int_{\R} x^{j_0} (\mu - \mu')(dx)$ $\neq$ $0$, otherwise, $\mu$ $=$ $\mu'$. In this case, define $p(\mu, x)$ $=$ $\left< x^{j_0}, \mu\right> $ $\in$ $\Hc_K$, then $p(\mu, x)$ $\neq$ $p(\mu', x)$. If $\mu$ $=$ $\mu'$, $x$ $\neq$ $x'$, take $p(\mu, x)$ $=$ $x$, then $p(\mu, x)$ $\neq$ $p(\mu', x')$. In either case, $\Hc_K$ separates points on $\Kc_K$ $\times$ $[-K, K]$.

\item {\it $\Hc_K$ vanishes at no point on $\Kc_K$ $\times$ $[-K, K]$}.
This is obvious as one can always choose a nonzero constant function in $\Hc_K$.
\end{itemize}

It thus follows from the Stone-Weierstrass lemma that $\Hc(\Kc_K \times [-K, K])$ is dense in $\Cc(\Kc_K \times [-K, K])$ with the topology of uniform convergence. Therefore, one can find a sequence of functions $p_n$, $ \tilde p_n$ $\in$ $\Hc_K$ such that for any $\epsilon$ $>$ $0$, there exists $N(\epsilon)$ such that for  $n$ $\geq$ $N(\epsilon)$
\beq \label{pnestim}
\sup_{(\mu, x) \in \Kc_K \times [-K, K]} \Big|p_n(\mu, x) - \partial_{xx} \frac{\delta \Phi}{\delta \mu}(\mu, x)\Big| \leq {\epsilon},\; \sup_{\mu \in \Kc_K} \Big|\tilde p_n(\mu) - \frac{\delta  \Phi}{\delta \mu}(\mu, 0)\Big| \leq \epsilon,
\enq

\medskip

\noindent {\it Step 2.} Let
\beqs
P_n(\mu, x) &:=& \tilde p_n(\mu)+ \int_{0}^x\int_{0}^y p_n(\mu, z) dzdy,\\
\Phi_n(\mu) &:=&  \Phi(\delta_0) + \int_0^1 \int_{[-K, K]} P_n(\lambda \mu + (1 -\lambda) \delta_0, x)(\mu - \delta_0) (dx) d\lambda.
\enqs
One can easily check that $\Phi_n$ $\in$ $\Gc(\Kc_K)$ by the above construction.
Now we have
\beqs
& & P_n(\mu, x) - \frac{\delta  \Phi}{\delta \mu}(\mu, x)\\
& & = \Big(\tilde p_n(\mu) + \int_{0}^x \int_{0}^y p_n(\mu, z) dzdy\Big) -  \Big(\frac{\delta  \Phi}{\delta \mu}(\mu, 0) + \int_{0}^x \int_{0}^y \partial_{xx} \frac{\delta  \Phi}{\delta \mu}(\mu, z)dzdy\Big)\\
&& = \tilde p_n(\mu)  -  \frac{\delta \Phi(\mu, 0)}{\delta \mu} +  \int_{0}^x \int_{0}^y \big(p_n(\mu, z) - \partial_{xx} \frac{\delta  \Phi}{\delta \mu}(\mu, z) \big)dzdy.
\enqs
Thus by \reff{pnestim},
\beqs
& & \sup_{\Kc_K \times [-K, K]}|\partial_x P_n(\mu, x) - \partial_\mu \Phi(\mu, x)| \leq K \cdot {\epsilon} ,  \\
& & \sup_{\Kc_K \times [-K, K]}|P_n(\mu, x) - \frac{\delta  \Phi}{\delta \mu}(\mu, x)| \leq (1 + K^2) \cdot {\epsilon}.
\enqs
Moreover,
\beqs
& & \Phi_n(\mu) - \Phi(\mu)\\
&& =  \Big(\Phi(\delta_0) + \int_0^1 \int_{[-K, K]}  P_n(\lambda \mu + (1 -\lambda) \delta_0, x)(\mu - \delta_0) (dx) d\lambda\Big)\\
& & \;\;- \Big(\Phi(\delta_0) + \int_0^1 \int_{[-K, K]} \frac{\delta  \Phi}{\delta \mu}(\lambda \mu + (1 -\lambda) \delta_0, x)(\mu - \delta_0) (dx) d\lambda\Big)\\
&& = \int_0^1 \int_{[-K, K]} \Big(P_n(\lambda \mu + (1 -\lambda) \delta_0, x) - \frac{\delta  \Phi}{\delta \mu}(\lambda \mu + (1 -\lambda) \delta_0, x) \Big)(\mu - \delta_0) (dx) d\lambda.
\enqs
Hence
\beqs
 \sup_{\Kc_K} \big|\Phi_n(\mu) - \Phi(\mu)\big| &\leq& 2(1 + K^2) \epsilon.
\enqs
Therefore,
\beqs
\|\Phi_n - \Phi\|_{\Kc_K} \leq  \epsilon \biggl(1 + K + 2(1 + K^2)\biggl),
\enqs
with $\Phi_n$ $\in$ $\Gc(\Kc_K)$, which is shown to be dense in $\Cc^{1, 1}(\Kc_K)$.
\ep

\vspace{2mm}

We are now ready to show that  It\^o's formula \reff{itoPP} holds for any $\Phi$ $\in$ $\Cc^{1, 1}(\Kc_K)$.
Without loss of generality,  assume that $X_0$ $=$ $0$. Fix $K >0$ and suppose $|X_r| \le K$ $\P$-a.s. for any $r \in [0,T]$. By Lemma \ref{Pdensity},
for any $\Phi \in \Cc^{1,1}(\mathcal{K}_K)$, there is a sequence of cylindrical functions $\{\Phi_n\}_{n=1}^\infty$  in $\Gc(\Kc_K)$ such that
\beq \label{PhinPhi}
\|\Phi_n - \Phi\|_{\mathcal{K}_K} \rightarrow 0.
\enq
By the definition of topology of $\Cc^{1, 1}$ in \reff{PhinPhi}, it is easy to deduce that
\begin{align} \label{Phinconv}
 \sup_{\mu \in \Kc_K} |\Phi_n(\mu) -\Phi(\mu)|  \to 0, \; & \sup_{\mu \in \Kc_K, X \sim \mu} \Big \lbrace |\partial_\mu \Phi_n(\mu, X) -\partial_\mu \Phi(\mu, X)|   \nonumber\\
& \hspace{2mm}  + |\partial_x \partial_\mu \Phi_n(\mu, X) -\partial_x\partial_\mu \Phi(\mu, X)|  \Big\rbrace \to 0 \ \ \ \P-a.s.
\end{align}
From \reff{lipsct1}, \reff{lipsct2}, and \reff{Phinconv}, there exists a  $K$-dependent  constant $C$ such that for all $n$ $\in$ $\N$,
\beq \label{derivbdd}
\sup_{\mu \in \Kc_K, X \sim \mu} \Big(, \partial_\mu \Phi_n(\mu, X)|, | \partial_\mu \Phi(\mu, X)|, |\partial_x\partial_\mu \Phi_n(\mu, X)|, |\partial_x\partial_\mu \Phi(\mu, X)|\Big) \leq C.
\enq

By Lemma \ref{M2C}, It\^o's formula holds  for $\{\Phi_n\}_{n \geq 1}$
\beq\label{itoPhin}
 & & \Phi_n(\P_{X_s }) - \Phi_n(\P_{X_t}) \nonumber\\
 & &={\mathbb{E}} \biggl[\int_t^s \partial_\mu \Phi_n(\P_{X_{r-} },{X}_{r-} ) d{X}_{r} + \frac{1}{2}\partial_x \partial_\mu \Phi_n (\P_{X_{r-} },{X}_{r-} ) d[{X} ,{X}]_r^c  \biggl] \nonumber\\
  & & \;+   \sum_{t < r \leq s} \Big\lbrace \Phi_n(\P_{X_{r}} ) - \Phi_n(\P_{X_{r-}} )\Big\rbrace 1_{\{\P_{X_r} \neq \P_{X_{r-}}\}} -   \E\Big[\sum_{t < r \leq s} \partial_\mu \Phi_n (\P_{X_{r-} },{X}_{r-} ) \Delta {X}_{r}\Big] \nonumber\\
  & & \; +  \E\Big[\sum_{t < r \leq s}\Big(\frac{\delta \Phi_n}{\delta \mu}(\P_{X_{r}}, X_{r}) - \frac{\delta \Phi_n}{\delta \mu} (\P_{X_{r}}, X_{r-})\Big) 1_{\{\P_{X_r} = \P_{X_{r-}}\}}\Big].
\enq
To establish a similar equation for any $\Phi$ $\in$ $\Cc^{1, 1}(\Kc_K)$, let us check that each term in the RHS of \reff{itoPhin} converges to its suitable limit.

\noindent $\bullet$ {\it First term of RHS of \reff{itoPhin}}.  Note from \reff{derivbdd} that the integrand is bounded. By applying the Dominated Convergence Theorem for stochastic integrals
\beqs
 \int_t^s \partial_\mu \Phi_n (\P_{X_{r-} },{X}_{r-} ) d{V}_{r} \overset{\text{a.s.}}{\longrightarrow} \int_t^s \partial_\mu \Phi(\P_{X_{r-} },{X}_{r-} ) d{V}_{r},\\
 \int_t^s \partial_\mu \Phi_n (\P_{X_{r-} },{X}_{r-} ) d{L}_{r} \overset{\text{p}}{\longrightarrow} \int_t^s \partial_\mu \Phi(\P_{X_{r-} },{X}_{r-} ) d{L}_{r},
\enqs
where $\overset{a.s}{\longrightarrow}$ and $\overset{\text{p}}{\longrightarrow}$ denote respectively the convergence $\P$-a.s. and in probability.\\
On one hand, from \reff{Xassum}
\beqs
\biggl|\int_t^s \partial_\mu \Phi_n (\P_{X_{r-} },{X}_{r-} ) d{V}_{r}\biggl| \leq C {\rm Var}(V)_s;
\enqs
and by \reff{Xassum} and the Dominated Convergence Theorem
\beq \label{integralV}
\lim_{n \to \infty}  \E \biggl[ \int_t^s \partial_\mu \Phi_n (\P_{X_{r-} },{X}_{r-} ) d{V}_{r}\biggl] =  \E\biggl[\int_t^s \partial_\mu \Phi(\P_{X_{r-} },{X}_{r-} ) d{V}_{r}\biggl].
\enq
On the other hand, by It\^o's isometry, \reff{Xassum}, and \reff{derivbdd}
\beqs
\sup_{n \geq 1}\E\biggl[ \biggl|\int_t^s \partial_\mu \Phi_n(\P_{X_{r-} },{X}_{r-} ) d{L}_{r}\biggl|^2\biggl] &=& \sup_{n \geq 1}\E\biggl[ \int_t^s \Big|\partial_\mu \Phi_n(\P_{X_{r-} },{X}_{r-} )\Big|^2 d [X, X]_r \biggl]\\
& \leq& C \E\big[[X, X]_s\big] < \infty.
\enqs
This implies the uniform integrability of $\big\lbrace\int_t^s (\partial_\mu \Phi_n )(\P_{X_{r-} },{X}_{r-} ) d{L}_{r}\big\rbrace_{n \geq 1}$ and
\beq \label{integralL}
\lim_{n \to \infty} \E \biggl[\int_t^s\partial_\mu \Phi_n(\P_{X_{r-} },{X}_{r-} ) d{L}_{r}\biggl] = \E\biggl[\int_t^s \partial_\mu \Phi(\P_{X_{r-} },{X}_{r-} ) d{L}_{r}\biggl].
\enq
Summing up \reff{integralV} and \reff{integralL}, we deduce that
\beqs
& & \lim_{n \to \infty} {\mathbb{E}} \biggl[\int_t^s \partial_\mu \Phi_n (\P_{X_{r-} },{X}_{r-} ) d{X}_{r}\biggl] = {\mathbb{E}} \biggl[\int_t^s \partial_\mu \Phi(\P_{X_{r-} },{X}_{r-} ) d{X}_{r}\biggl].
\enqs
\noindent $\bullet$ {\it Second term of RHS of \reff{itoPhin}}.  By \reff{derivbdd},  $\big|\int_t^s \partial_x \partial_\mu \Phi_n (\P_{X_{r-} },{X}_{r-} ) d[{X} ,{X}]_r^c\big|$ $\leq$ $C [X, X]_s^c$. By the Dominated Convergence Theorem,
\beqs
& & \lim_{n \to \infty} {\mathbb{E}}\biggl[\int_t^s \partial_x \partial_\mu \Phi_n (\P_{X_{r-} },{X}_{r-} ) d[{X} ,{X}]_r^c \biggl]={\mathbb{E}}\biggl[\int_t^s \partial_x \partial_\mu \Phi (\P_{X_{r-} },{X}_{r-} ) d[{X} ,{X}]_r^c \biggl].\\
\enqs
\noindent $\bullet$ {\it Third term of RHS of \reff{itoPhin}}.  Note that $|\Delta X_r|$ $\leq$ $2K$ $\P$-a.s.. By \reff{Phinconv},
\beqs
& & \lim_{n \to \infty}\Big\lbrace\Big(\Phi_n(\P_{X_r}) - \Phi_n(\P_{X_{r-}})\Big)1_{\{\P_{X_r} \neq \P_{X_{r-}}\}} - \partial_\mu \Phi_n(\P_{X_{r-}},  X_{r-}) \Delta  X_r\Big\rbrace\\
& & = \Big(\Phi(\P_{X_r}) - \Phi(\P_{X_{r-}})\Big)1_{\{\P_{X_r} \neq \P_{X_{r-}}\}} - \partial_\mu \Phi(\P_{X_{r-}},  X_{r-}) \Delta  X_r, \;\;\; \P-a.s..
\enqs
Moreover, by Remark \ref{Itoremainder} and \reff{derivbdd}
\beqs
& & \sup_{n \geq 1}\E\biggl[\sum_{t < r\leq s}\Big\lbrace\Big(\Phi_n(\P_{X_r}) - \Phi_n(\P_{X_{r-}})\Big)1_{\{\P_{X_r} \neq \P_{X_{r-}}\}} - \partial_\mu \Phi_n(\P_{X_{r-}},  X_{r-}) \Delta  X_r\Big\rbrace\biggl]\\
& &  \leq 2 \sup_{n \geq 1} \sup_{(\mu, x) \in \Hc_{K}}|\partial_\mu \Phi_n(\mu, x)|  \E\Big[\sum_{t < r \leq s}|\Delta  X_r|\Big] \leq  C \E\Big[\sum_{t < r \leq s}|\Delta  X_r|\Big].
\enqs
Applying the Dominated Convergence Theorem yields,
\beqs
& & \lim_{n \to \infty} \E\biggl[\sum_{t < r \leq s} \Big\lbrace\Big(\Phi_n(\P_{X_r}) - \Phi_n(\P_{X_{r-}})\Big)1_{\{\P_{X_r} \neq \P_{X_{r-}}\}} - \partial_\mu \Phi_n(\P_{X_{r-}},  X_{r-}) \Delta  X_r\Big\rbrace\biggl]\\
& & = \E\biggl[\sum_{t < r \leq s}\Big\lbrace\Big(\Phi(\P_{X_r}) - \Phi(\P_{X_{r-}})\Big)1_{\{\P_{X_r} \neq \P_{X_{r-}}\}} - \partial_\mu \Phi(\P_{X_{r-}},  X_{r-}) \Delta  X_r \Big\rbrace \biggl].
\enqs
\noindent $\bullet$ {\it Fourth term of RHS of \reff{itoPhin}}. By \reff{derivbdd}, $$\sup_n \E\Big[\sum_{t < r \leq s}\Big|\Big(\frac{\delta \Phi_n}{\delta \mu}(\P_{X_{r}}, X_{r}) - \frac{\delta \Phi_n}{\delta \mu} (\P_{X_{r}}, X_{r-})\Big)\Big| 1_{\{\P_{X_r} = \P_{X_{r-}}\}}\Big] < \infty.$$
By the Dominated Convergence Theorem,
\beqs
& & \lim_{n \to \infty} \E\Big[\sum_{t < r \leq s}\Big|\Big(\frac{\delta \Phi_n}{\delta \mu}(\P_{X_{r}}, X_{r}) - \frac{\delta \Phi_n}{\delta \mu} (\P_{X_{r}}, X_{r-})\Big)\Big| 1_{\{\P_{X_r} = \P_{X_{r-}}\}}\Big]\\
&  & =\E\Big[\sum_{t < r \leq s}\Big|\Big(\frac{\delta \Phi}{\delta \mu}(\P_{X_{r}}, X_{r}) - \frac{\delta \Phi}{\delta \mu} (\P_{X_{r}}, X_{r-})\Big)\Big| 1_{\{\P_{X_r} = \P_{X_{r-}}\}}\Big].
\enqs
Now taking the limit on both sides of \reff{itoPhin}, we establish It\^o's formula \reff{itoPP} for any $\Phi$ $\in$ $\Cc^{1, 1}(\Kc_K)$.

\subsubsection{Proof of Theorem \ref{TheoremIto}}

We now finish the proof of Theorem \ref{TheoremIto}, using a localization argument  for general functions over $\Pc_2(\R)$.

Take a general process ${X}$ $=$ $(X_t)_{t \in [0, T]}$, and note that ${X}_{-}$ $:= {(X_{t-})}_{t \in [0, T]}$ is locally bounded.
Fix $K>0$, denote $\tau_K$ $=$ $\inf\{ 0 < t \leq T$, $|X_{t-}| > K\}$, and define the truncated process
\beqs
X_t^{K} = X^{\tau_K}_{t-} = X_{t-  \wedge \tau_K}.
\enqs
Then $|X_t^{K}|$ $\leq$ $K$ $\P$-a.s.. Now, for any $\Phi$ $\in$ $\Cc^{1, 1}(\Pc_2(\R))$, one can apply It\^o's formula to $\Phi|_{\Kc_K}$, the restriction of $\Phi$ over $\Kc_K$:
  \beq \label{ItoPhiKk}
& & \Phi(\P_{X^{K}_s }) - \Phi(\P_{X^{K}_t})\\
& & =
   {\mathbb{E}} \biggl[ \int_t^s \partial_\mu \Phi (\P_{X^{K}_{r-} },{X}^{K}_{r-} ) d{X}^{K}_{r} + \frac{1}{2} \int_t^s \partial_{x} \partial_\mu \Phi(\P_{X^{K}_{r-} },{X}^{K}_{r-} ) d[{X}^{K} ,{X}^{K}]_r^c\biggl] \nonumber \\
& &+\; \sum_{t < r \leq s}  \Big\lbrace  \Phi(\P_{X^{K}_{r}} ) - \Phi(\P_{X^{K}_{r-}} )\Big\rbrace  1_{\{\P_{X^K_r} \neq \P_{X^K_{r-}}\}} -   {\mathbb{E}} \Big[\sum_{t < r \leq s}\partial_\mu \Phi (\P_{X^{K}_{r-} },{X}^{K}_{r-} ) \Delta {X}^{K}_{r}\Big] \nonumber\\
& & + \;  \E\Big[\sum_{t < r \leq s}\Big(\frac{\delta \Phi}{\delta \mu}(\P_{X^K_{r}}, X^K_{r}) - \frac{\delta \Phi}{\delta \mu} (\P_{X^K_{r}}, X^K_{r-})\Big) 1_{\{\P_{X^K_r} = \P_{X^K_{r-}}\}}\Big].  \; \nonumber
\enq
To see that It\^o's formula \reff{itoPP} holds for any $\Phi$ $\in$ $\Cc^{1, 1}(\Pc_2(\R))$, let us check each term of the RHS of \reff{ItoPhiKk}.

\noindent $\bullet$ {\it First term of RHS of \reff{ItoPhiKk}}. By the construction of $X^{K}$ and the stopping rule for stochastic integral
\beqs
& & \int_t^s \partial_\mu \Phi(\P_{X^{K}_{r-} },{X}^{K}_{r-} ) d{X}^{K}_{r} = \int_t^s \partial_\mu \Phi(\P_{X_{r-}},{X}_{r-} )  1_{[0, \tau_K]} d{X}_{r},
\enqs
Since $|(\partial_\mu \Phi )(\P_{X_{r-}},{X}_{r-} )  1_{[0, \tau_K]}|$ $\leq$ $|(\partial_\mu \Phi )(\P_{X_{r-}},{X}_{r-} )|$ $\leq$ $C$, $\P$-a.s., by the Dominated Convergence Theorem for stochastic integrals (Proposition 2.74 \cite{M2007})
\beqs
\int_t^s \partial_\mu \Phi(\P_{X_{r-}},{X}_{r-} )  1_{[0, \tau_K]} d{V}_{r} \overset{a.s.}{\longrightarrow} \int_t^s \partial_\mu \Phi(\P_{X_{r-} },{X}_{r-} ) d{V}_{r},\\
\int_t^s \partial_\mu \Phi(\P_{X_{r-}},{X}_{r-} )  1_{[0, \tau_K]} d{L}_{r} \overset{p}{\longrightarrow} \int_t^s \partial_\mu \Phi(\P_{X_{r-} },{X}_{r-} ) d{L}_{r}.
\enqs
Repeating the same argument as the proof for the first term of RHS of \reff{itoPhin},
\beq \label{step2first}
& & \lim_{K \to \infty} {\mathbb{E}} \biggl[ \int_t^s \partial_\mu \Phi(\P_{X^{K}_{r-} },{X}^{K}_{r-} ) d{X}^{K}_{r}\biggl] =  {\mathbb{E}} \biggl [ \int_t^s \partial_\mu \Phi(\P_{X_{r-} },{X}_{r-} ) d{X}_{r} \biggl].
\enq
\noindent $\bullet$ {\it Second term of RHS of \reff{ItoPhiKk}}.  Note that we have in the pathwise sense
\beqs
& & \int_t^s \partial_{x} \partial_\mu \Phi(\P_{X^{K}_{r-} },{X}^{K}_{r-} ) d[{X}^{K} ,{X}^{K}]_r^c =\int_t^s \partial_x\partial_\mu \Phi (\P_{X_{r-} },{X}_{r-} ) 1_{[0, \tau_K]}d[{X} ,{X}]_r^c\\
& & \overset{a.s.}{\longrightarrow} \int_t^s \partial_x\partial_\mu \Phi (\P_{X_{r-} },{X}_{r-} )d[{X} ,{X}]_r^c.
\enqs
Since $|\partial_x\partial_\mu \Phi(\mu, x)|$ $\leq$ $C$ and $\E[[X, X]_T^c]$ $<$ $\infty$ for any $T$ $>$ $0$, Dominated Convergence Theorem implies
\beqs \label{step2second}
& & \lim_{K \to \infty} \E \biggl[\int_t^s \partial_x \partial_\mu \Phi(\P_{X^{K}_{r-} },{X}^{K}_{r-} ) d[{X}^{K} ,{X}^{K}]_r^c\biggl]  = {\mathbb{E}}\biggl [ \int_t^s \partial_x\partial_\mu \Phi(\P_{X_{r-} },{X}_{r-} ) d[{X} ,{X}]_r^c\biggl].
\enqs
\noindent $\bullet$ {\it Third term of RHS of \reff{ItoPhiKk}}.
By the construction of $X_t^K$,
\beqs
 & & \Big(\Phi(\P_{X^{K}_{r}} ) - \Phi(\P_{X^{K}_{r-}})\Big) 1_{\{\P_{X^K_r} \neq \P_{X^K_{r-}}\}} -  \partial_\mu \Phi(\P_{X^{K}_{r-} },{X}^{K}_{r-} ) \Delta {X}^{K}_{r} \\
 &  & \overset{a.s.}{\longrightarrow} \Big(\Phi(\P_{X_{r}} ) - \Phi(\P_{X_{r-}} )\Big)  1_{\{\P_{X_r} \neq \P_{X_{r-}}\}}-  \partial_\mu \Phi(\P_{X_{r-} },{X}_{r-} ) \Delta {X}_{r}.
\enqs
By Remark \ref{Itoremainder},  RHS of the above equation is integrable under $\E[\sum_{t < r \leq s}\cdot]$, and by the Dominated Convergence Theorem
\beqs  \label{step2third}
& & \lim_{K \to \infty }  {\mathbb{E}} \Big[\sum_{t < r \leq s}  \Big\lbrace  \Big(\Phi(\P_{X^{K}_{r}} ) - \Phi(\P_{X^{K}_{r-}})\Big) 1_{\{\P_{X^K_r} \neq \P_{X^K_{r-}}\}} -  \partial_\mu \Phi (\P_{X^{K}_{r-} },{X}^{K}_{r-} ) \Delta {X}^{K}_{r} \Big \rbrace \Big] \nonumber\\
& & =  {\mathbb{E}} \Big[ \sum_{t < r \leq s} \Big\lbrace  \Big(\Phi(\P_{X_{r}} ) - \Phi(\P_{X_{r-}})\Big) 1_{\{\P_{X_r} \neq \P_{X_{r-}}\}}-  \partial_\mu \Phi (\P_{X_{r-} },{X}_{r-} ) \Delta {X}_{r}\Big\rbrace\Big].
\enqs
\noindent $\bullet$ {\it Fourth term of RHS of \reff{ItoPhiKk}}. As
\beqs
\Big(\frac{\delta \Phi}{\delta \mu}(\P_{X^K_{r}}, X^K_{r}) - \frac{\delta \Phi}{\delta \mu} (\P_{X^K_{r}}, X^K_{r-})\Big)1_{\{\P_{X^K_r} = \P_{X^K_{r-}}\}} \overset{a.s.}{\longrightarrow} \Big(\frac{\delta \Phi}{\delta \mu}(\P_{X_r}, X_{r}) - \frac{\delta \Phi}{\delta \mu} (\P_{X_{r}}, X_{r-}) \Big)1_{\{\P_{X_r} = \P_{X_{r-}}\}},
\enqs
by the Dominated Convergence Theorem,
\beq \label{step2four}
& & \lim_{K \to \infty} \E\Big[\sum_{t < r \leq s} \Big(\frac{\delta \Phi}{\delta \mu}(\P_{X^K_{r}}, X^K_{r}) - \frac{\delta \Phi}{\delta \mu} (\P_{X^K_{r}}, X^K_{r-})\Big)  1_{\{\P_{X^K_r} \neq \P_{X^K_{r-}}\}}\Big] \nonumber\\
& & = \E\Big[\sum_{t < r \leq s} \Big(\frac{\delta \Phi}{\delta \mu}(\P_{X_r}, X_{r}) - \frac{\delta \Phi}{\delta \mu} (\P_{X_{r}}, X_{r-})\Big) 1_{\{\P_{X_r} = \P_{X_{r-}}\}}\Big].
\enq
Now, taking limit on both sides of \reff{ItoPhiKk} and $\Phi(\P_{X_t^K})$ $\to$ $\Phi(\P_{X_t})$ as $K$ $\to$ $\infty$, It\^o's formula \reff{itoPP} holds for any $\Phi$ $\in$ $\Cc^{1, 1}(\Pc_2(\R))$.
\ep


\subsection{Proof of Theorem \ref{Theoremito_diffcond} }

Note that all the arguments of Theorem \ref{TheoremIto} remain essentially the  same here. The only modification is to ensure the Dominated Convergence Theorem holds whenever appropriate.  Specifically, instead of relying on $|\partial_\mu \Phi(\mu, x)| \leq C$ and $|\partial_x\partial_\mu \Phi(\mu, x)| \leq C$ in Definition \ref{defregularityf} for \reff{step2first}-\reff{step2four} in the proof of Theorem \ref{TheoremIto}, one first needs to check that  conditions \reff{conditionPhi} in Definition \ref{defregularityf_relax} and \reff{strongerX} in  {\bf (H)$_{strict}$} are sufficient for the Dominated Convergence Theorem. That is to check
\beqs
& &\E\biggl[\Big|\int_t^s  \partial_\mu \Phi(\P_{X_{r-}},{X}_{r-} ) d V_r\Big|\biggl]^2
\leq \sup_{0 < r \leq T} \E\Big[|\partial_\mu \Phi (\P_{X_{r-}},{X}_{r-} )|^2\Big] \cdot \E\big[|{\rm Var}(V)_T|^2\big]\\
& &\E\biggl[\Big|\int_t^s \partial_\mu \Phi (\P_{X_{r-} },{X}_{r-} ) d[X, X]_r \Big|\biggl]^2
\leq\sup_{0 < r \leq T} \E\Big[|\partial_\mu \Phi (\P_{X_{r-}},{X}_{r-} )|^2\Big] \cdot \E\big[|[X, X]_s|^2\big],\\
& &  \E\Big[\sum_{t < r \leq s} \Big|\Big(\Phi(\P_{X_r}) - \Phi(\P_{X_{r-}})\Big)1_{\{\P_{X_r} \neq \P_{X_{r-}}\}} - \partial_\mu\Phi(\P_{X_{r-}},  X_{r-}).\Delta X_r \Big|\Big] ^2\\
& \leq  &  4\sup_{0 < r \leq T} \E\Big[|\partial_\mu \Phi (\P_{X_{r-}},{X}_{r-} )|^2\Big] \E\Big[\Big(\sum_{0 < t \leq T} |\Delta X_t|\Big)^2\Big],
\enqs
and
\beqs
& & \E\Big[\sum_{t < r \leq s} \Big(\frac{\delta \Phi}{\delta \mu}(\P_{X_r}, X_{r}) - \frac{\delta \Phi}{\delta \mu} (\P_{X_{r}}, X_{r-})\Big)1_{\{\P_{X_r} = \P_{X_{r-}}\}}\Big]^2  \\
& \leq &  \sup_{0 < r \leq T} \E\Big[|\partial_\mu \Phi (\P_{X_{r-}},{X}_{r-} )|^2\Big] \cdot \E\Big[\Big(\sum_{0 < t \leq T} |\Delta X_r|\Big)^2\Big],
\enqs
where \reff{Phiintegralform}-\reff{equ: Itolinearlions} in Remark \ref{Itoremainder} are used. Now, it remains to verify conditions \reff{strongerX}-\reff{conditionPhi} are satisfied: given a semimartingale $X$, $\widehat\Kc: = \Kc_X \cup \Kc_{X-}$ is compact in $\Pc_2(\R)$, with $\Kc_X = \big\{\P_{X_r}, 0 \leq r \leq T\big\}$ and $\Kc_{X-} = \big\{\P_{X_{r-}}, 0 < r \leq T\big\}$. To see this, let $\{\mu_{r_n}\}_{n \geq 1}$ be any sequence in $\widehat \Kc$, where $\mu_{r_n}$ is either $\P_{X_{r_n}} \in \Kc_X$ or $\P_{X_{r_{n-}}} \in \Kc_{X-}$. Since $\{r_n\}_{n \geq 1}$ has a convergent subsequence in $[0, T]$, denoted as $\{r_{n_k}\}_{k \geq 1}$,  one can show that $\mu_{r_{n_k}}$ also converges with its limit in $\widehat\Kc$. Consequently by \reff{conditionPhi}, $\sup_{0 < r \leq T} \E\Big[|\partial_\mu \Phi (\P_{X_{r-}},{X}_{r-} )|^2\Big] \leq \sup_{\widehat \Kc} \E\Big[|\partial_\mu \Phi (\P_{X_{r-}},{X}_{r-} )|^2\Big]$ $< \infty$.
\ep



\begin{Remark}[Generalization to $d > 1$] \label{remd1}
For ease of exposition, the proof for Theorem \ref{TheoremIto} is given for $d=1$. Nevertheless, its adaptation for the case of $d > 1$ is straightforward, including  Definition \ref{defGM}, Definition \ref{defHMU}, Lemma \ref{M2C}, and Lemma \ref{Pdensity}. Indeed by \cite{CLSF2019}, the underlying space in $\Gc(\Mc)$ in Definition \ref{defGM} and $\Hc(\Mc, \Uc)$ in Definition \ref{defHMU} can be any locally compact Polish space.  When the underlying space is $\R^d$, the compact set in $\Pc_2(\R^d)$ can be defined as
$\Kc_K: = \{\mu  \in \Pc_2(\R^d)| \mathrm{supp}(\mu) \subset [-K, K]^d\}$. One can also check that the proof of Lemma \ref{M2C} remains unchanged except for the polynomial $g_k$, $1 \leq k \leq n$, now  defined from $\R^d$ to $\R$; hence the derivative of $g_k$ is replaced by its partial derivative with respect to $x_i$, $1 \leq i \leq d$.

Finally, for Lemma \ref{Pdensity}, one can find ${\bf p_n}(\mu, x): = (p_n^{i, j}(\mu, x))_{1 \leq i, j \leq d}$ with $p_n^{i, j}(\mu, x) = p_n^{j, i}(\mu, x)$ such that $\big|p_n^{i, j} (\mu, x) - \partial_{x_i}\partial_{x_j} \frac{\delta \Phi}{\delta \mu} \Phi(\mu, x)\big| \leq \epsilon$, and then obtain $P_n(\mu, x)$ by integrating ${\bf p}_n(\mu, x)$ with respect to $dz = dz^1\ldots dz^d$ and $dy = dy^1\ldots dy^d$. The rest proceeds similarly as the proof of Lemma \ref{Pdensity}.
 \end{Remark}

\begin{Remark} \label{comaretaltouzha21} 
Note that
we propose two different sets of conditions for It\^o's formula: one in Theorem \ref{TheoremIto} with  stronger conditions on the functional and  weaker integrability conditions on the semimartingale, the other in Theorem \ref{Theoremito_diffcond} with weaker conditions on the functional and stronger integrability condition on the semimartingale.
Note that  conditions in   \cite{taltouzha21} are  similar to ours in Theorem \ref{Theoremito_diffcond}, as both sets of conditions allow for   the interchange of the expectation and integration.

Comparing conditions \reff{conditionPhi} and \reff{strongerX} with Definition 3.1 and equation (3.2) in \cite{taltouzha21}, it is clear  that our condition \reff{strongerX} on the semimartingale $X$ is from the Doob-Meyer decomposition of the semimartingale $X$ while equation (3.2) in \cite{taltouzha21} is from the Doob-Meyer decomposition of the continuous part of the semimartingale $X$.
Note that we assume  $\partial_\mu \Phi$ and $\partial_x\partial_\mu \Phi$ in \reff{conditionPhi} (Definition \ref{defregularityf_relax}) to be bounded in $L^2$, locally uniform in $\mu$, whereas \cite{taltouzha21} assumes that (see their Definition 3.1) $\frac{\delta \Phi}{\delta \mu}$ has quadratic growth in $x$, $\partial_\mu \Phi$ has linear growth in $x$, and $\partial_x\partial_\mu \Phi$ is bounded, all locally uniform in $\mu$; The assumptions in Definition 3.1,  \cite{taltouzha21} provide with sufficient conditions for our assumptions in \reff{conditionPhi}. 
\end{Remark}

\begin{Remark}\label{remark:particle}
 We would also like to point out that there appears to have some essential difficulties to adopt  the particle approximation approach for general semimartingales. For instance,
 there are two terms for jumps in It\^o's formula in Theorems \ref{TheoremIto} \& \ref{Theoremito_diffcond}: the first jump term sums over the set of time points at which the law of the semimartingale is not continuous, and the second jump term sums over the set of time points at which the law of the semimartingale is continuous while the semimartingale itself has jumps. Using the particle approximation approach, the first jump term in the standard It\^o's formula of the empirical projection converges to the first jump term, as expected. However, it is unclear whether the second jump term in the standard It\^o's formula of the empirical projection converges to the second jump term associated with the linear derivative in Theorem \ref{TheoremIto} or Theorem \ref{Theoremito_diffcond}.
 To see this more precisely,  suppose that $((X_t^n)_{t \geq 0})_{n \geq 1}$ is a sequence of i.i.d. copies of the semimartingale $(X_t)_{t \geq 0}$. Denote by
$\mu_t^N = \frac{1}{N} \sum_{n=1}^N \delta_{X_t^n}$ as the empirical measure of  $((X_t^n)_{t \geq 0})_{n \geq 1}$. Given an integer $N \geq 1$, the empirical projection of $\Phi:$ $\mathcal{P}_2(\mathbb{R}^d) \to \R$ onto $\R^d$ is defined as
\beqs
\phi^N: (\R^d)^N \ni (x_1, \ldots, x_N) \mapsto \Phi(\frac{1}{N} \sum_{i=1}^N \delta_{x_i}).
\enqs
The strategy of the particle approximation approach is then to expand $\phi^N: (\R^d)^N \ni (x_1, \ldots, x_N)$ using classical It\^o's formula and take the limit.  However, one cannot see clearly in the expansion of $\phi^N$ whether the term $\E[\sum_{t < r \leq s}\Big(\phi^N(X_r^1, \ldots, X_r^N) - \phi^N(X_{r-}^1, \ldots, X_{r-}^N)\Big) 1_{\{\P_{X_r} = \P_{X_{r-}}\}}]$  converges, as $N \to \infty$, to $\E[\sum_{t < r \leq s}\Big(\frac{\delta \Phi}{\delta \mu}(\P_{X_{r}}, X_{r}) - \frac{\delta \Phi}{\delta \mu} (\P_{X_{r}}, X_{r-}) \Big)  1_{\{\P_{X_r} = \P_{X_{r-}}\}} ]$ in \reff{itoPP}.
\end{Remark}


\vspace{1mm}

%

The remaining part of this paper demonstrates  how   It\^o's formula enables us to derive  dynamic programming equations and    verification theorems
 for  McKean--Vlasov controls with jump diffusions and for McKean--Vlasov mixed regular-singular control problems.
 It also allows for generalizing the classical relation  between the maximum principle and the dynamic programming principle  to the  McKean--Vlasov singular control setting, where the adjoint process is expressed in terms of the derivative of the value function with respect to probability measures.

\section{McKean--Vlasov Control with Jump-diffusion Process} \label{secMKVjump}

\subsection{Problem Formulation and HJB Equation} \label{subsecMKVjump}

On a given probability space $(\Omega, \Fc, \F, \P)$, we consider a $d$-dimensional Brownian motion $W = (W_t)_{t \geq 0}$, and an independent Poisson random measure $N$ with a finite intensity measure $\nu$, {with} $\tilde N(d\theta, dt): = N(d\theta, dt) - \nu(d\theta)dt$ the compensated Poisson random measure. Fix $T < \infty$, consider the following McKean--Vlasov control problem with jump-diffusion process, where the  $\R^d$-valued state variable $X_s$ ($s \geq t$) follows the dynamics starting from $t \in [0, T]$,
\begin{equation} \label{equMKVSDE}
\left\{
\begin{array}{rcl}
d X_s &=& b(X_s, \alpha_s, \P_{X_{s}^{}})ds + \sigma(X_s, \alpha_s, \P_{X_s}) d W_s \\
& & \;\;\; \hspace{2cm} + \int_{\R^q} \beta(X_{s-}, \alpha_{s^-}, \P_{X_{s-}}, \theta) N(ds, d\theta), \; s \geq t, \\
X_{t-} &=& \xi \in L^2(\Fc_t; \R^d),
\end{array}
\right.
\end{equation}
Here the control $\alpha$ $=$ $(\alpha_s)_{s \geq t}$ is a c\`adl\`ag  $\F$-adapted process taking values in a subset $A$ of $\R^m$, satisfying the square integrability  condition:
$\E\bigl[\int_t^T|\alpha_s|^2ds\bigl] < \infty$, and denoted by $\alpha$ $\in$ $\Ac_t$.

Define a cost functional
\beq \label{equMKVcost}
J(t, \xi, \alpha) = \E\biggl[\int_t^T f(X_s^{t,\xi,\alpha}, \alpha_s, \P_{X_s^{t,\xi,\alpha}})ds + g(X_T^{t,\xi,\alpha}, \P_{X_T^{t,\xi,\alpha}})\biggl],
\enq
then the McKean--Vlasov control problem is to find an optimal control (if exists)  to minimize the cost functional $J$. That is to solve for
\beq \label{equMKVvalue}
V(t, \xi) = \inf_{\alpha \in \Ac_t} J(t, \xi, \alpha),
\enq
subject to the jump diffusion \reff{equMKVSDE}.

To ensure the well-definedness of the control problem, the following conditions on coefficients $b$, $\sigma$ and $\beta$ are imposed.

\noindent {\bf (H1)} There exists a constant $C_{b, \sigma, \beta} >0$ such that the coefficients $b$, $\sigma$ defined from $\R^d \times A \times \Pc_2(\R^d)$ to $\R^d$, $\R^{d \times d}$ respectively and $\beta$ from $\R^d \times {A} \times \Pc_2(\R^d) \times \R^q$ to $\R^d$ satisfy:
\beqs
|b(x, a, \mu)|+ |\sigma(x, a, \mu)| + \int_{\R^q}|\beta(x, a, \mu, \theta)|\nu(d\theta) \leq C_{b, \sigma, \beta}\Big(1 + |x| + |a| + \|\mu\|_2\Big),
\enqs
for any $x \in \R^d$, $a \in A$, and $\mu \in \Pc_2(\R^d)$.
Moreover, there exists a constant $C'_{b, \sigma, \beta}>0$ such that
\beqs
& & |b(x, a, \mu) - b(x', a', \mu')| + |\sigma(x, a, \mu) - \sigma(x', a', \mu')| + \int_{\R^q}|\beta(x, a, \mu, \theta) - \beta(x', a', \mu', \theta)|\nu(d\theta)\\
 &\leq& C'_{b, \sigma, \beta} \Big(|x-x'| + |a-a'| + W_2(\mu, \mu')\Big),
\enqs
 for any $x, x'$ $\in$ $\R^d$, $a, a' \in A$, and $\mu, \mu' \in \Pc_2(\R^d)$.

One can easily check that if $\alpha \in \Ac_t$, the dynamics \reff{equMKVSDE} under conditions {\bf (H1)} has a unique square integrable solution, and we denote by $X$ the solution to \eqref{equMKVSDE} and omit the superscripts of $X^{t, \xi, \alpha}$ whenever there is no confusion.

Furthermore, to ensure that the cost functional in \reff{equMKVcost} is well defined and finite, we will make the following assumptions on coefficients $f$ and $g$.

\medskip

\noindent {\bf (H2)}\; There exists $C_{f, g} > 0$ such that the coefficients $f$ $: \R^d \times A \times \Pc_2(\R^d) \to \R$ and $g: \R^d \times \Pc_2(\R^d) \to \R$ satisfy:
\beqs
|f(x, a, \mu)| + |g(x, \mu)| \leq C_{f, g} \Big(1 + |x|^2 + |a|^2 + \|\mu\|_2^2 \Big),
\enqs
 for all $x$ $\in$ $\R^d$, $a \in A$, and $\mu \in \Pc_2(\R^d)$.

\medskip

Now,  proceeding by the same argument as in \cite{cospha19}, we have the dynamic programming principle (PP) for this McKean--Vlasov control problem \reff{equMKVSDE}-\reff{equMKVvalue}.

\vspace{2mm}

\paragraph{DPP}  Under assumption {\bf (H1)}, the value function in \reff{equMKVvalue} is law-invariant.
Moreover, for $\mu$ $=$ $\P_\xi$ $\in$ $\Pc_2(\R^d)$, with a slight abuse of notation, write
\beq  \label{equMKVvalue_law_invariant}
 V(t,\mu) = V(t,\xi) = \inf_{\alpha \in \Ac_t} J(t, \xi, \alpha).
\enq
Then, the dynamic programming principle holds such that
\begin{align} \label{DPPjump}
V(t, \mu) = \inf_{\alpha \in \Ac_t} \Big \lbrace \E\Big[\int_t^s f(X_r^{t, \xi, \alpha}, \alpha_r, \P_{X_r^{t, \xi, \alpha}})dr \Big] + V(s, \P_{X_s^{t, \xi, \alpha}})\Big\rbrace, \tag{MKV-Jump-DPP}
\end{align}
for all $t, s \in [0, T]$, with $t < s$, $\mu \in \Pc_2(\R^d)$ and any $\xi \in L^2(\Fc_t; \R^d)$ such that $\P_{\xi} = \mu$.

\vspace{2mm}

Next, define an operator $H$  on $\R^d  \times A \times \Pc_2(\R^d)  \times \R^d \times \R^{d \times d} \times \L_\nu$, with $\L_\nu$  the set of all deterministic integrable functions with respect to $\nu$, such that
\beq \label{equ_jump_H}
H(x, a, \mu, p, M, r(\cdot)) &=&f(x, a, \mu) + b(x, a, \mu).p + \frac{1}{2}{\rm Tr}\big(\sigma\sigma\trans(x, a, \mu) M\big)\\
& & \;\; + \;\int_{\R^q} \big(r(x + \beta(x, a, \mu, \theta)) -r(x)\big)\nu(d\theta). \nonumber
\enq
Then, if the value function in \reff{equMKVvalue} is sufficiently smooth, applying Corollary \ref{CorItojump} to \reff{DPPjump} yields the following HJB equation
\beq \label{equjumpHJB}
\left\{
\begin{array}{rcl}
\partial_t V(t, \mu) + \E\Big[\Inf_{a \in A} H(\xi, a, \mu, \partial_\mu V(t, \mu, \xi), \partial_x\partial_\mu V(t, \mu, \xi), \frac{\delta  V}{\delta \mu}(t, \mu, \xi))\Big]&=&0,\\
V(T, \mu) &=& \langle g(\cdot, \mu), \mu \rangle.
\end{array}
\right.
\enq


One can further connects the value function $V(t, \mu)$ in \reff{DPPjump} with the HJB equation \reff{equjumpHJB} in the following verification theorem.

\subsection{Verification Theorem}
\begin{Theorem} [Verification Theorem]\label{Propjumpver}
Let $\widehat V: [0, T] \times \Pc_2(\R^d) \to \R$ be a function in $\Cc^{1, (1, 1)}([0, T] \times \Pc_2(\R^d))$.
\begin{itemize}
\item [(a)] Suppose that $\widehat V$ is a solution to the HJB equation \reff{equjumpHJB}, then $\widehat V(t, \mu) \leq V(t, \mu)$ for any $(t, \mu) \in [0, T] \times \Pc_2(\R^d)$, with $V$ given in \reff{equMKVvalue_law_invariant}.
\item [(b)] Moreover, suppose that  there exists a continuous function $\hat a(t, x, \mu)$ for $(t, x, \mu)$ $\in$  $[0, T)$ $\times$ $\R^d$ $\times$ $\Pc_2(\R^d)$, valued in $A$ such that
\beq \label{hata}
\hat a (t, x, \mu) \in \arg\min_{a \in A} H(x, a, \mu, \partial_\mu \widehat V(t, \mu, x), \partial_x\partial_\mu \widehat V(t, \mu, x), \frac{\delta  \widehat V}{\delta \mu}(t, \mu, x));
\enq
suppose also that the corresponding controlled dynamics under \reff{hata}
\beqs
\left\{
\begin{array}{rcl}
dX_s &=&  b(X_s,\hat a(s, X_s, \P_{X_s}),  \P_{X_s}) ds + \sigma(X_s, \hat a(s, X_s, \P_{X_s}), \P_{X_s})dW_s \\
& & \hspace{0.4cm} + \int_{\R^q} \beta(X_{s-}, \hat a(s, X_{s^-}, \P_{X_{s}}), \P_{X_{s}}, \theta) N(ds, d\theta), \; t \leq s \leq T, \\
X_{t-} &=& \xi.
\end{array}
\right.
\enqs
admits a unique solution denoted as $X^*_s$,
 and the control $\alpha^*$ $=$ $\{\alpha^*_s: = \hat a (s, X_s^*, \P_{X_s^*}), t \leq s < T\} \in \Ac_t$. Then $\widehat V =V$,  with $V$ given in \reff{equMKVvalue_law_invariant}, and $\alpha^*$ is an optimal Markovian control.
\end{itemize}
\end{Theorem}
{\bf Proof.} (a) Fix $(t, \xi) \in [0, T] \times L^2(\Fc_t; \R^d)$, with $\P_{\xi} =\mu \in \Pc_2(\R^d)$, and consider $X$ $=$ $X^{t, \xi, \alpha}$, a solution to SDE \eqref{equMKVSDE}
for an arbitrary control $\alpha \in \Ac_t$. By assumption ({\bf H1}),
\beqs
\E\Big[\sup_{t \leq s \leq T} |X_s^{t, \xi, \alpha}|^2\Big] \leq C_T \big(1 + \E[|\xi|^2]\big) < \infty,
\enqs
which implies that assumption ({\bf H}) holds, i.e.,
\beqs
& & \E\biggl[\int_t^T \biggl\lbrace|b(X_s,  \alpha_s, \P_{X_s})|^2 + |\sigma(X_s, \alpha_s, \P_{X_s})|^2 +  \int_{\R^d}|\beta(X_s, \alpha_s, \P_{X_s}, \theta)| \nu(d\theta) \biggl\rbrace ds\biggl] < \infty.
\enqs
Now applying Corollary \ref{CorItojump} to $\hat V(s, \P_{X_s^{t, \xi, \alpha}})$ between $t$ and $T$ yields
\beqs
 \widehat V(t, \mu) &=& \widehat V(T, \P_{X_T})-\int_t^T \biggl\lbrace \partial_s \widehat V(s, \P_{X_s}) + \E\Big[\partial_\mu \widehat V(\P_{X_s}, X_s). b(X_s, \alpha_s, \P_{X_s})\\
  & & \;\;+\; \frac{1}{2} {\rm Tr}\big(\sigma\sigma\trans(X_s, \alpha_s, \P_{X_s})\partial_x\partial_\mu \widehat V(\P_{X_s}, X_s)\big)\\
  & & \;\; + \; \int_{\R^q} \Big( \frac{ \delta \widehat V}{\delta \mu}(\P_{X_{s}}, X_{s} + \beta(X_s, \alpha_s, \P_{X_s}, \theta)) - \frac{\delta \widehat V}{\delta \mu}(\P_{X_{s}}, X_{s})\Big)\nu(d\theta)\Big] \biggl\rbrace ds\\
   &=& \widehat V(T, \P_{X_T}) -\int_t^T \biggl\lbrace \partial_s \widehat V(s, \P_{X_s}) - \E\Big[f (X_s, \alpha_s, \P_{X_s})\Big]\\
   & & \;\;+\; \E\Big[H(X_s, \alpha_s, \P_{X_s}, \partial_\mu \widehat V(s, \P_{X_s}, X_s), \partial_x\partial_\mu \widehat V(s, \P_{X_s}, X_s), \frac{\delta \widehat V}{\delta \mu}(s, \P_{X_s}, X_s))\Big]\biggl \rbrace ds\\
 & \leq & \widehat V(T, \P_{X_T^{}}) + \E\bigg[\int_t^T f (X_s^{}, \alpha_s, \P_{X_s^{}}) ds\bigg] \\
 &=& \E\bigg[ g(X_T,\P_{X_T}) + \int_t^T f (X_s^{}, \alpha_s, \P_{X_s^{}}) ds\bigg]
 \\
 & \leq & V(t, \mu),
\enqs
where the first inequality is from the HJB equation \reff{equjumpHJB} for
$\widehat V$,  the last equality is from the terminal condition for $\widehat V$, and the final inequality is by the definition of $V$.

\vspace{1mm}

\noindent (b) Now we apply the same argument again with the control $\alpha^*$ $\in$ $\Ac_t$  given by  $\alpha^*_s = \hat a(s, X_s^*, \P_{X_s^*})$, $t\leq s\leq T$, with $\hat a$ attaining the infimum in \reff{hata}, and obtain
\beqs
 \widehat V(t, \mu) = \widehat V(T, \P_{X_T^{t, \xi, \alpha^*}}) + \E\biggl[\int_t^T f (X_s^{t, \xi, \alpha^*}, \alpha^*_s, \P_{X_s^{t, \xi, \alpha^*}})\biggl]dr,
\enqs
which implies that $\widehat V(t, \mu) = V(t, \mu)$.
\ep

\begin{Remark}\label{remark4.4}
In the case when there is no smooth solution to the HJB equation \reff{equjumpHJB}, then its solution should be interpreted in the viscosity sense as in
\cite{BIVS2019}, \cite{CGKPR20}, \cite{CGKPR21}.  Notice that
a particular form of  It\^o's formula for cylindrical functions has also used in
\cite{BIVS2019} in order to derive the viscosity solution property of the value function.
\end{Remark}

\subsection{Example: McKean--Vlasov Linear Quadratic Control Problem with Jump-diffusion Process} \label{secMKVLQ}

In the particular case of a McKean--Vlasov linear quadratic (LQ) control problem with jump diffusion, we will illustrate how to apply the Bellman equation \reff{equjumpHJB} and the verification theorem \ref{Propjumpver} to derive explicit solutions.
In fact, a similar problem has been solved  in \cite{TM2019} via a duality approach where the jump is driven by a Poisson random martingale measure, and also studied in \cite{BDT2019} which considers jump-diffusion-regime switching state dynamics of multiple decisions makers.

For simplicity, take $d=1$ and $A=\R$ as the analysis can be easily generalized to the multivariate case. In the LQ case, coefficients of dynamics in \reff{equMKVSDE} are specified as
\beqs
b(x, a, \mu) &=& b_0 + b_1 x + \bar b_1 \bar \mu + b_2 a,\\
\sigma(x, a, \mu)&=& \sigma_0 + \sigma_1 x + \bar \sigma_1 \bar \mu + \sigma_2 a,\\
\beta(x, a, \mu, \theta) &=& \beta_0(\theta) + \beta_1(\theta)x + \bar\beta_1(\theta) \bar\mu + \beta_2(\theta) a,
\enqs
for $(x, a, \mu) \in \R \times \R \times \Pc_2(\R)$, where $\bar\mu: = \int_{\R} x \mu(dx)$;
the cost functional in \reff{equMKVcost} takes the form of
\beqs
f(x, a, \mu) &=& f_1 x^2 + \bar f_1 \bar \mu^2 + f_2 a^2,\\
g(x, \mu) &=&  g_1 x^2 + \bar g_1 \bar\mu^2,
\enqs
where ${\rm Var}(\mu): = \int x^2 \mu(dx) - \bar\mu^2$. Here $b_0, \sigma_0, b_i, \sigma_i, (i=1, 2), \bar b_1, \bar\sigma_1$,
$f_1, \bar f_1, f_2, g_1, \bar g_1$, are constants in $\R$; and $\beta_0, \beta_1, \bar\beta_1, \beta_2$ are deterministic functions of $\theta$.

We start by guessing the form of the solution to the HJB equation \reff{equjumpHJB} with
\beq \label{LQexampleW}
\widehat V(t,\mu) = A(t){\rm Var}(\mu) + B(t) \bar\mu^2 + C(t) \bar\mu + D(t)
\enq
for some time-dependent deterministic functions $A(t), B(t), C(t)$ and $D(t)$.
Then
\beqs
\frac{\delta \widehat V}{\delta \mu} (t, \mu, x) &=& A(t) x^2 + 2 \Big(B(t) - A(t)\Big) \bar\mu x +  C(t) x  + c,\\
\partial_\mu \widehat V(t, \mu, x) &=& 2 A(t) x + 2 \Big(B(t) - A(t)\Big) \bar\mu +C(t),\\
\partial_x\partial_\mu \widehat V(t, \mu, x) &=& 2 A(t),
\enqs
where $c$ is a constant.
Suppose that $\widehat V$ in \reff{LQexampleW} satisfies the HJB equation \reff{equjumpHJB}, by some straightforward calculations, one can decompose $H$ in \reff{equ_jump_H} into two parts, with one part involving $a$, denoted as $G(a)$, and the other part independent of $a$, such that
\beqs
&  &\inf_{a \in \R}H(x, a, \mu,  \partial_\mu \widehat V(t, \mu), \partial_x \partial_\mu \widehat V(t, \mu), \frac{ \delta \widehat V}{\delta \mu}(t, \mu))\\
&=& \inf_{a \in \R} G(a) + \Big(f_1 + 2 b_1 A(t) +  |\sigma_1|^2 A(t)  + \langle 2\beta_1 + |\beta_1|^2, \nu \rangle A(t) \Big) x^2\\
& & + \Big(2 (B(t) - A(t)) b_1 + 2 \bar b_1 A(t)+ 2\sigma_1\bar\sigma_1 A(t) + 2 \langle \beta_1 \bar\beta_1 + \bar\beta_1, \nu \rangle A(t) + 2 \langle \beta_1, \nu\rangle (B(t) -A(t)) \Big) x \bar \mu \\
& & \;+\;  \Big(\bar f_1 +  2 \bar b_1(B(t) - A(t))  + |\bar\sigma_1|^2 A(t)  +  \langle |\bar\beta_1|^2, \nu \rangle A(t) + 2\langle \bar\beta_1, \nu \rangle (B(t) -A(t))\Big)\bar\mu^2 \\
& & \;+ \; \Big(2 b_0 A(t)  + b_1  C(t) + 2\sigma_0\sigma_1 A(t) + 2\langle\beta_0\beta_1 + \beta_0, \nu\rangle A(t) + \langle \beta_1, \nu\rangle C(t) \Big) x\\
& & \;+\; \Big(\bar b_1 C(t) + 2b_0(B(t) - A(t)) + 2\sigma_0\bar\sigma_1A(t)  + 2\langle \beta_0 \bar\beta_1, \nu\rangle A(t) + 2\langle \beta_0, \nu \rangle (B(t) -A(t))\\
&  &\;+\; \langle \bar\beta_1, \nu\rangle C(t)\Big)\bar\mu +  b_0 C(t) + \sigma_0^2 A(t) + \langle |\beta_0|^2, \nu\rangle A(t) + \langle \beta_0, \nu\rangle C(t),
\enqs
where $G(a)$ takes the following quadratic form
\beq \label{equGa}
G(a) &=& U(t) a^2  + 2S(t) x a  + 2(Z(t) - S(t))\bar\mu a + Y(t) a\\
&=& U(t) \Big( a + \frac{S(t)}{U(t)} (x -\bar\mu) + \frac{Z(t)}{U(t)} \bar\mu + \frac{Y(t)}{2U(t)}\Big)^2 \nonumber\\
& & \;-\; \frac{S(t)^2}{U(t)} (x  -\bar\mu)^2 - \frac{Z(t)^2}{U(t)} \bar\mu^2 -\frac{Z(t)Y(t)}{U(t)} \bar\mu- \frac{Y(t)^2}{4U(t)}   \nonumber\\
&  &\; -\; \frac{2 S(t) Z(t)}{U(t)} (x - \bar\mu) \bar\mu  - \frac{S(t)Y(t)}{U(t)} (x-\bar\mu)  \nonumber,
\enq
with
\beqs
\left\{
\begin{array}{rcl}
U(t)  &=& f_2  + \sigma_2^2 A(t) + \langle \beta_2^2, \nu \rangle A(t), \\
S(t) &= & b_2 A(t) + \sigma_1\sigma_2 A(t) + \langle \beta_1\beta_2 + \beta_2, \nu\rangle A(t),\\ 
Z(t) &=&  b_2 B(t) + (\sigma_1 + \bar\sigma_1) \sigma_2 A(t) + \langle (\beta_1 + \bar\beta_1)\beta_2, \nu \rangle A(t) + \langle \beta_2, \nu \rangle B(t),\\
Y(t) &= & b_2 C(t) + 2\sigma_0 \sigma_2 A(t) + 2 \langle \beta_0\beta_2, \nu \rangle A(t) + \langle \beta_2, \nu\rangle  C(t).\\ 
\end{array}
\right.
\enqs
It is clear by \reff{equGa} that $G(a)$ attains the infimum at
\beq \label{LQa*}
\hat a(t, x, \mu)= -\frac{S(t)}{U(t)} (x -\bar\mu) - \frac{Z(t)}{U(t)} \bar\mu - \frac{Y(t)}{2U(t)}.
\enq
Substituting \reff{equGa} with $\hat a$ into \reff{equjumpHJB}, we obtain
\beq \label{LQHJBequation}
 & & 0= \partial_t \widehat V(t, \mu) + \E\Big[\inf_{a \in \R} H(\xi, a, \mu,  \partial_\mu \widehat V(t, \mu, \xi), \partial_x \partial_\mu \widehat V(t, \mu, \xi), \frac{ \delta \widehat V }{\delta \mu}(t, \mu, \xi))\Big] \\
& & =\Big(\dot A(t) + f_1 + 2 b_1 A(t) + |\sigma_1|^2  A(t) + \langle 2\beta_1 + |\beta_1|^2, \nu \rangle A(t) - \frac{S(t)^2}{U(t)}\Big) {\rm Var}(\mu) \nonumber\\
& & \;\;\;+\; \Big(\dot B(t) + f_1 + \bar f_1 + 2(b_1 + \bar b_1) B(t) + (\sigma_1 + \bar \sigma_1)^2 A(t) + \langle (\beta_1 + \bar\beta_1)^2, \nu \rangle A(t) \nonumber\\
&  & \hspace{3cm}+\; 2\langle (\beta_1 + \bar\beta_1), \nu\rangle B(t)  -\frac{Z(t)^2}{U(t)}\Big)\bar\mu^2 \nonumber\\
&  & \;\;\; +\;  \Big(\dot C(t) + 2 (b_1 + \bar b_1) C(t) + 2 b_0 B(t) + 2 \sigma_0(\sigma_1 + \bar \sigma_1)A(t) + 2 \langle \beta_0(\beta_1 + \bar\beta_1), \nu \rangle A(t) \nonumber\\
& &  \hspace{3cm} + \; 2 \langle \beta_0, \nu\rangle B(t) + \langle \beta_1 + \bar\beta_1, \nu\rangle C(t) \Big)\bar\mu \nonumber\\
& & \;\;\; + \; \dot D(t) +  b_0 C(t) + \sigma_0^2 A(t) + \langle |\beta_0|^2, \nu \rangle A(t) -\frac{Y(t)^2}{4 U(t)}.  \nonumber
\enq

Now comparing terms in ${\rm Var}(\mu)$, $\bar\mu^2$, $\bar\mu$ in \reff{LQHJBequation}, we get the following ODEs system for $A(t)$, $B(t)$, $C(t)$ and $D(t)$,
\beqs
\left\{
\begin{array}{rcl}
\dot A(t) + f_1 + 2 b_1 A(t) + |\sigma_1|^2  A(t) + \langle 2\beta_1 + |\beta_1|^2, \nu \rangle A(t) - \frac{S(t)}{U(t)^2} =0,\\ A(T) = g_1,
\end{array}
\right.
\enqs
\beqs
\left\{
\begin{array}{rcl}
\dot B(t) + f_1 + \bar f_1 + 2(b_1 + \bar b_1) B(t) + (\sigma_1 + \bar \sigma_1)^2 A(t) + \langle (\beta_1 + \bar\beta_1)^2, \nu \rangle A(t) \\
+ 2\langle (\beta_1 + \bar\beta_1), \nu\rangle B(t)  -\frac{Z(t)^2}{U(t)} =0,\\
B(T) = g_1 + \bar g_1,
\end{array}
\right.
\enqs
\beqs
\left\{
\begin{array}{rcl}
\dot C(t) + 2 (b_1 + \bar b_1) C(t) + 2 b_0 B(t) + 2 \sigma_0(\sigma_1 + \bar \sigma_1)A(t) + 2 \langle \beta_0(\beta_1 + \bar\beta_1), \nu \rangle A(t) \nonumber\\
+ 2 \langle \beta_0, \nu\rangle B(t) + \langle \beta_1 + \bar\beta_1, \nu\rangle C(t) =0, \\
C(T) = 0,
\end{array}
\right.
\enqs
\beqs
\left\{
\begin{array}{rcl}
\dot D(t) + 2 b_0 C(t) + \sigma_0^2 A(t) +\langle |\beta_0|^2, \nu\rangle A(t) + \langle \beta_0, \nu\rangle C(t) -\frac{Y(t)^2}{4U(t)} =0,\\
D(T) =0.
\end{array}
\right.
\enqs

By the Verification Theorem \ref{Propjumpver}, the solution for the  McKean--Vlasov LQ control problem is given by the solution of the Riccati equations for $A(t)$ and $B(t)$ and those of the linear equations for $C(t)$ and $D(t)$ given $A(t)$ and $B(t)$. See \cite{Yong2013} for sufficient conditions on the existence of solutions to Riccati equations.

From \reff{LQa*}, the optimal control is given by
\beqs
\alpha^*_t =  -\frac{S(t)}{U(t)} (X_t^* - \E[X_t^*]) - \frac{Z(t)}{U(t)} \E[X_t^*] - \frac{Y(t)}{2U(t)},
\enqs
where $X^*_t$ is the controlled dynamic associated with $\alpha^*_t$. In fact, one can further compute $\E[X_t^*]$ so that
\beqs
d\E[X_t^*] = R(t)\E[X_t^*]dt + Q(t),
\enqs
with
\beqs
R(t) &=&  b_1 + \bar b_1 + \langle \beta_1 + \bar\beta_1, \nu\rangle  - \frac{Z(t)}{U(t)} (b_2 + \langle \beta_2, \nu \rangle,\\
Q(t) &=& b_0 + \langle \beta_0, \nu \rangle - \frac{Y(t)}{2U(t)} (b_2 + \langle \beta_2, \nu\rangle).
\enqs
Note that the path of $X_t^*$ is discontinuous, yet $\E[X_t^*]$ is continuous.

\section{McKean--Vlasov Mixed Regular-Singular Control Problem} \label{secMKVsingular}

\subsection{Problem Formulation and Dynamic Programming Equation}

On a given probability space $(\Omega, \Fc, \F, \P)$, we consider a $d$-dimensional Brownian motion $W = (W_t)_{t \geq 0}$. A mixed regular and singular controlled McKean--Vlasov dynamics can be formulated as the following SDE starting from $t \in [0, T]$
\begin{align}\label{equdynamicssingular}
\left\{
\begin{array}{rcl}
d{X_s}  &= &{b} (X_s, \alpha_s, \P_{{X_{s}}} ) ds + {\sigma}(X_s,  \alpha_s, \P_{{X_{s}}})  d{W_s} +  {\lambda} d\eta_s, \; t \leq s \leq T,\\
 X_{t-} &=& \xi \in L^2(\Fc_t; \R^d),
 \end{array}
 \right.
\end{align}
for some measurable functions $b: \R^d \times \R^m \times \Pc_2(\R^d) \to \R^d$, $\sigma: \R^d \times \R^m \times \Pc_2(\R^d) \to \R^{d \times d}$, and some nonnegative constant $\lambda ={\rm diag}(\lambda_1, \ldots, \lambda_d) \in \R^{d \times d}$.
Here  the mixed regular and singular control pair $(\alpha_s, \eta_s)_{t \leq s \leq T}$ is a pair  of  processes valued in $\R^m$ $\times [0, \infty)^d$ such that
\begin{itemize}
\item  $(\alpha_s)_{t \leq s \leq T} \in \Ac_t:= \Big\{\alpha: [t, T] \to \R^m  \; \big| \; \alpha \text{ is } \F-\text{adapted}, \text{ and } \E\biggl[\int_t^T|\alpha_s|^2ds\biggl] < \infty\Big\}$,
\item $(\eta_s)_{t \leq s \leq T}$ $\in \Uc_t$ is a nondecreasing c\`adl\`ag $\F$-adapted process satisfying $\E[{\rm Var}(V)_T] < \infty$.
\end{itemize}

Given the dynamics \reff{equdynamicssingular},  the  McKean--Vlasov regular-singular control problem is to derive,  over the admissible control set $\Ac_t$ $\times$ $\Uc_t$, the following value function
\begin{align}\label{MKVsingularcost}
\begin{split}
\tilde{V}(t, \xi)  & =  \inf_{(\alpha, \eta)\in \Ac_t \times \Uc_t} \tilde{J}(t, \xi, \alpha, \eta)
\\  &= \inf_{(\alpha, \eta) \in \Ac_t \times \Uc_t} \mathbb{E} \biggl [\int_t^T \Big( f(X_s, \alpha_s, \P_{X_{s}} ) ds  + g(X_T, \P_{X_T})\Big)  + \int_{[t, T)} \gamma d \eta_s \biggl],
\end{split}
\end{align}
where $\gamma$ is a fixed constant in $\R^d$ and $\int_{[t, T)} \gamma d \eta_t$ means that a jump at the terminal time $T$ is not allowed.

To ensure this McKean--Vlasov control problem \reff{equdynamicssingular}-\reff{MKVsingularcost} is well defined, we will make the following assumptions on the coefficients for the coefficients $b$ and $\sigma:$

\medskip

\noindent {\bf (H3)} There exists a constant $C_{b, \sigma} > 0$ such that
 for any $x, x' \in \R^d$, $a, a' \in \R^m$, and $\mu, \mu' \in \Pc_2(\R^d)$
\beqs
|b(x, a, \mu) - b(x', a', \mu')| + |\sigma(x, a, \mu) - \sigma(x', a', \mu')| &\leq& C_{b, \sigma} \big(|x - x'| + W_2(\mu, \mu') + |a-a'|\big),
\enqs
and
\beqs
|b(x, a, \mu)| + |\sigma(x, a, \mu)| \leq C_{b, \sigma}\big(1 + |x| + \|\mu\|_2 + |a|\big).
\enqs


 Note that under condition ({\bf H3}), equation \reff{equdynamicssingular} has a unique strong solution, denoted as $X^{t, \xi, \alpha, \eta}$,  for any admissible control $(\alpha, \eta) \in \Ac_t \times \Uc_t$.

Furthermore, we will assume that the running functions $f$ and $g$ satisfy the square growth condition.

\noindent {\bf (H4)} There exists a constant $C_{f, g} > 0$ such that
\beqs
|f(x, a, \mu)| + |g(x, \mu)| \leq C_{f, g} \Big(1 + |x|^2 + |a|^2 + \|\mu\|_2^2\Big),
\enqs
for all $x$ $\in$ $\R^d$, $a \in \R^m$, $\mu \in \Pc_2(\R^d)$


\vspace{3mm}

Next, proceeding by the same arguments as \cite{CGKPR20}, we can obtain the DPP for this McKean--Vlasov  regular-singular control problem \reff{equdynamicssingular}-\reff{MKVsingularcost}.

\vspace{3mm}

Under assumptions {\bf (H3)}-{\bf (H4)}, for any $t, s$ $\in$ $[0, T]$ with $t < s$, and $\xi \in L^2(\Fc_t; \R^d)$ with $\P_{\xi} =\mu$, we have
\begin{align} \label{equDPP}
V(t, \mu) \triangleq \tilde V(t, \xi)\;=\;  \inf_{(\alpha, \eta) \in \Ac_t \times \Uc_t} \E\Big[\int_t^{s} f(X_r, \alpha_r, \P_{X_{r}}) dr + \int_{[t, s]}\gamma d \eta_r\Big]  + V(s, \P_{X_s^{t, \xi}}). \tag{MKV-Regular-Singular-DPP}
\end{align}

\medskip


Suppose that $V(t, \mu)$ is sufficiently smooth. By applying Corollary \ref{CorItosingular} to \reff{equDPP}, we obtain the dynamic programming equation:
\beq \label{equsingularDP}
0  &=& \inf_{(\alpha, \eta) \in \Ac_t \times \Uc_t} \int_t^s \biggl\lbrace \partial_r V(r, \P_{X_r}) + \E\Big[\partial_\mu V(r, \P_{X_r}, X_r).b(X_r, \alpha_r, \P_{X_r})\\
& & \; +\; \frac{1}{2}{\rm Tr}\big(\sigma\sigma\trans(X_r, \alpha_r, \P_{X_r})\partial_x\partial_\mu V(r, \P_{X_r}, X_r)\big) + f(X_r, \alpha_r, \P_{X_r})\Big]\biggl\rbrace dr \nonumber\\
&  &\;+ \;\E\Big[\int_t^s  \big(\lambda \partial_\mu V(r, \P_{X_{r-}}, X_{r-}) + \gamma\big) d\eta_r\Big] - \E\Big[\sum_{t \leq r \leq s} \partial_\mu V(r, \P_{X_{r-}}, X_{r-}). \Delta X_r\Big] \nonumber \\
& & \;+ \;\sum_{t \leq r \leq s} \Big(V(r, \P_{X_r}) - V(r, \P_{X_{r-}})\Big)1_{\{\P_{X_r} \neq \P_{X_{r-}}\}}  \nonumber\\
&  & \; + \;  \E\Big[\sum_{t \leq r \leq s}\Big(\frac{\delta V}{\delta \mu}(r, \P_{X_{r}}, X_{r}) - \frac{\delta V}{\delta \mu} (r, \P_{X_{r}}, X_{r-})\Big) 1_{\{\P_{X_r} = \P_{X_{r-}}\}}\Big] \nonumber,
\enq
 for all $0 \leq t < s < T$. Notice that the last two terms of RHS of \reff{equsingularDP} can be written in integral form
\beqs \label{singluar_jump_integral}
V(r, \P_{X_r}) - V(r, \P_{X_{r-}}) &=& \int_0^1 \E[\partial_\mu V(r, \P_{X_{r-} + h\Delta X_r}, X_{r-} + h\Delta X_{r}).\Delta X_r] dh,\\\label{singular_jump2_integral}
\frac{\delta V}{\delta \mu}(r, \P_{X_{r}}, X_{r}) - \frac{\delta V}{\delta \mu} (r, \P_{X_{r-}}, X_{r-}) &=& \int_0^1 \partial_\mu V(r, \P_{X_{r-} + h \Delta X_r}, X_{r-} + h \Delta X_r).\Delta X_r dh.
\enqs
Hence, by rearranging the last three terms in \reff{equsingularDP}, we  obtain the following definition of a classical solution to \reff{equsingularDP}.

\begin{Definition} Consider $\widehat V \in \Cc^{1, (1, 1)}([0, T] \times \Pc_2(\R^d))$, define
\beq \label{equnonaction}
\Cc(\widehat V): = \Big\{(t, \mu, x) \in [0, T] \times \Pc_2(\R^d) \times \R^d: \gamma + \lambda \partial_\mu \widehat V(t, \mu, x) > 0\Big\},
\enq
denote by $\Cc_{t, \mu}(\widehat V)$ the projection of $\Cc(\widehat V)$ onto $(t, \mu)$-space, and define
\beq \label{operatorHsingular}
H(x, a, \mu, p, M) \;=\;f(x, a, \mu) + b(x, a, \mu).p + {\rm Tr}\big(\sigma(x, a, \mu)\trans M \big).
\enq
We say that $\widehat V$ is a classical solution to the dynamic programming equation \reff{equsingularDP} if the following conditions hold
\begin{itemize}
\item for every $(t, \mu) \in [0, T] \times \Pc_2(\R^d)$, if $(t, \mu, \xi) \in \Cc(\widehat V)$, $\P$-a.s. with $\xi \sim \mu$.
\beq \label{equHJBnonaction}
\partial_t \widehat V(t, \mu) + \E\big[ \inf_{a \in \R^m}  H(\xi, a, \mu, \partial_\mu \widehat V(t, \mu, \xi),\partial_x\partial_\mu \widehat V(t, \mu, \xi) \sigma(\mu, a, \xi))\big]=0,
\enq
\item for every $\mu \in \Pc_2(\R^d), \widehat V(T, \mu) = \langle g(\cdot, \mu), \mu \rangle$,
\item for every $(t, \mu) \in [0, T] \times \Pc_2(\R^d)$ and $\xi \sim \mu$,
\beq \label{equpartialW}
 \gamma + \lambda \partial_\mu \widehat V(t, \mu, \xi) \geq 0,\;\;\; \P-a.s.,
\enq
\item for every $(t, \mu, a)$ $\in [0, T] \times \Pc_2(\R^d) \times \R^m$,
\beq \label{equWvariational}
\partial_t \widehat V(t, \mu) + \E \big[ H(\xi, a, \mu, \partial_\mu \widehat V(t, \mu, \xi), \partial_x\partial_\mu \widehat V(t, \mu, \xi)\sigma(\mu, a, \xi))\big]\geq 0.
\enq
\end{itemize}
\end{Definition}

\subsection{Verification Theorem}
\begin{Theorem} [Verification Theorem] \label{veri} \;\;
\begin{itemize}
\item [(a)]
Suppose that $\widehat V$ $ \in \mathcal{C}^{1, (1,1)}([0, T] \times \mathcal{P}_2(\mathbb{R}^d))$ is a classical solution of the dynamic programming equation \reff{equsingularDP},
then $\widehat V(t, \mu) \leq V(t, \mu)$ for any $(t, \mu)$ $\in$ $[0, T] \times \Pc_2(\R^d)$, where $V$ is given in \reff{equDPP}.
\item [(b)] Furthermore, if there exists $(\alpha^*, \eta^*) \in \Ac_t \times \Uc_t$ such that
\beq \label{equX*nonaction}
\P\Big(\text {Leb a.e. } s\in [t, T), X_s^*: = X_s^{t, \xi, \alpha^*, \eta^*} \in \Cc_{s, \P_{X_s^*}}(\widehat V)\Big) = 1,
\enq
and for every $(s, \P_{X_s^*}, X_s^*) \in \Cc(\widehat V)$,
\beq \label{equsingularalpha*}
\alpha^*_s = \hat a(s, X_s^*, \P_{X_s^*}), \text{ with}\;  \hat a(s, x, \mu)\in \arg\min_{a \in \R^m} H(x, a, \mu, \partial_\mu \widehat V(s, \mu), \partial_x\partial_\mu \widehat V(s, \mu)),
\enq
\beq \label{equeta*}
\E \int_{[t, T)} \Big(\lambda \partial_\mu \widehat V(s, \P_{X^*_{s-}}, X_{s-}^*) + \gamma\Big) d \eta^*_s =0,
\enq
\beq \label{DPPWX*}
 \widehat V(s, \P_{X_{s-}^{*}})  = \widehat V(t, \P_{X_s^{*}})  + \gamma \mathbb{E}  \big[\Delta \eta_s^* \big] \;\;\; \mbox{for any} \; s \in \big\{ t \leq s < T: \P_{X_s^*} \neq \P_{X_{s-}^*}\big\},
\enq
and
\begin{align}\label{equlinearhatVX*}
\P\Big(\frac{\delta \widehat V}{\delta \mu}(s, \P_{X_{s}^*}, X_{s}^*) = \frac{\delta \widehat V}{\delta \mu} (s, \P_{X_{s}^*}, X_{s-}^*)  + \gamma \Delta \eta_s^*, \; &\mbox{for all}\;t \leq s < T \nonumber \\
 &\mbox{with}\; \P_{X_s^*} = \P_{X_{s-}^*} \Big) = 1.
\end{align}
Then $\widehat V$ is the value function of problem \reff{equdynamicssingular}-\reff{MKVsingularcost}. That is, $\widehat V(t, \mu) = V(t, \mu)$ for any $(t, \mu) \in [0, T] \times \Pc_2(\R^d)$, where $V$ is given in \reff{equDPP}.
\end{itemize}
\end{Theorem}
{\bf Proof}.\; (a) Fix $(t,\xi) \in [0, T] \times L^2(\Fc_t; \R^d)$ with $\P_{\xi} =\mu$ and consider $X_s$ $=$ $X_s^{t, \xi, \alpha, \eta}$ solution to SDE \reff{equdynamicssingular} for a given arbitrary admissible control $(\alpha, \eta) \in \Ac_t \times \Uc_t$. Under conditions ({\bf H3})-{\bf (H4)},
\beqs
\E\Big[\sup_{t \leq s \leq T} |X_s^{t, \xi, \alpha, \eta}|^2\Big] \leq C_T \big(1 + \E[|\xi|^2]\big) < \infty.
\enqs
This implies assumption ({\bf H}).
Applying Corollary \ref{CorItosingular} to $\widehat V(s, \P_{X_s^{t, \xi}})$ between $t$ and $T$,
\beq
\widehat V(t, \mu) &=& \widehat V(T, \P_{X_T}) - \int_t^T \biggl\lbrace \partial_s \widehat V(s, \P_{X_s}) +  \E\Big[\partial_\mu \widehat V(s, \P_{X_s}, X_s).b(X_s, \alpha_s, \P_{X_s}) \nonumber\\
& & \; + \; \frac{1}{2} {\rm Tr}\big(\sigma\sigma\trans(X_s, \alpha_s, \P_{X_s})\partial_x\partial_\mu \widehat V(s, \P_{X_{s}}, X_{s})\big)\biggl]ds \biggl\rbrace
\label{proofequItoW} \\
& & \; - \; \E\Big[\int_t^T \lambda \partial_\mu \widehat V(s, \P_{X_{s-}}, X_{s-}). d\eta_s^c \Big] \nonumber\\
& & \;-\; \sum_{t \leq s < T} \Big(\widehat V(s, \P_{X_s}) - \widehat V(s, \P_{X_{s-}})\Big) 1_{\{\P_{X_s} \neq \P_{X_{s-}}\}}\nonumber\\
& & \;-\;  \E\Big[\sum_{t \leq s < T} \Big(\frac{\delta \widehat V}{\delta \mu}(s, \P_{X_{s}}, X_{s}) - \frac{\delta \widehat V}{\delta \mu} (s, \P_{X_{s}}, X_{s-})\Big) 1_{\{\P_{X_s} = \P_{X_{s-}}\}}   \Big] \nonumber.
\enq
Since $\widehat V$ is a classical solution of \reff{equsingularDP}, by \reff{equWvariational},
\beqs
& & \partial_s \widehat V(s, \P_{X_s}) +  \E\Big[\partial_\mu \widehat V(s, \P_{X_s}, X_s).b(X_s, \alpha_s, \P_{X_s}) +  \frac{1}{2} {\rm Tr}\big(\sigma\sigma\trans(X_s, \alpha_s, \P_{X_s})\partial_x\partial_\mu \widehat V(s, \P_{X_{s}}, X_{s})\big)\Big]\\
& & \geq - f(X_s, \alpha_s, \P_{X_{s}}).
\enqs
Moreover, by \reff{equpartialW}, we have
\beqs
\E\Big[\int_t^T \lambda \partial_\mu \widehat V(s, \P_{X_{s-}}, X_{s-}). d\eta_s^c \Big]  \geq - \gamma E\Big[\int_t^T d\eta_s^c\Big].
\enqs
According to the Fundamental Theorem of Calculus and \reff{equpartialW},
\beqs
& & \widehat V(s, \P_{X_s}) - \widehat V(s, \P_{X_{s-}}) = \int_0^1 \E \Big[\lambda \partial_\mu \widehat V(s, \P_{X_{s-} + h\Delta X_s}, X_{s-} + h\Delta X_s) . \Delta  \eta_s\Big]dh \geq -\gamma \E[\Delta \eta_s],\\
& &  \frac{\delta \widehat V}{\delta \mu}(r, \P_{X_{r}}, X_{r}) - \frac{\delta \widehat V}{\delta \mu} (r, \P_{X_{r}}, X_{r-}) = \int_0^1 \partial_\mu \widehat V(r, \P_{X_r}, X_{r-} + h \Delta X_r).\Delta X_r dh \geq -\gamma \Delta \eta_s.
\enqs
From the above three inequalities, with the terminal condition for $\widehat V$ plugged into \eqref{proofequItoW}, we obtain
\beqs
\widehat V(t, \mu) \leq \E\biggl[\int_t^T \Big( f(X_s, \alpha_s, \P_{X_{s}} ) dt  + g(X_T, \P_{X_T})\Big) + \int_{[t, T)} \gamma d \eta_s \biggl].
\enqs
\medskip

\noindent (b) Applying \reff{proofequItoW} with $\P_{X_s^*}$, we see from \reff{equnonaction}, \reff{equX*nonaction}, and \reff{equsingularalpha*} that the second term on the right hand side of \reff{proofequItoW} becomes
\beqs
\E\Big[\int_t^T f(X_s^*, \alpha^*_s, \P_{X_s^*})ds\Big].
\enqs
Furthermore, by \reff{equeta*}-\reff{equlinearhatVX*},
\beqs
& & - \E\Big[\int_t^T \lambda \partial_\mu \widehat V(s, \P_{X_{s-}^*}, X_{s-}^*) d\eta_s^{*,c} \Big] -\sum_{t \leq s < T} \Big(\widehat V(s, \P_{X_s^*}) - \widehat V(s, \P_{X_{s-}^*})\Big) 1_{\{\P_{X_s}^* \neq \P_{X_{s-}^*}\}}\\
& & -  \E\Big[\sum_{t \leq s < T} \Big(\frac{\delta \widehat V}{\delta \mu}(s, \P_{X_{s}^*}, X_{s}^*) - \frac{\delta \widehat V}{\delta \mu} (s, \P_{X_{s}^*}, X_{s-}^*)\Big) 1_{\{\P_{X_s^*} = \P_{X_{s-}^*}\}}  \Big] \\
& = &  \E\Big[\int_t^T \gamma  d\eta_s^{*,c}\Big] + \gamma \sum_{t \leq s < T} \E[\Delta \eta_s^*]1_{\{\P_{X_s^*} \neq \P_{X_{s-}^*}\}} +  \gamma \E\Big[\sum_{t \leq s < T} \Delta \eta_s^* 1_{\{\P_{X_s^*} = \P_{X_{s-}^*}\}}\Big] \\
& = &  \E\Big[\int_t^T \gamma d \eta_s^*\Big].
\enqs
Now by applying \reff{proofequItoW} to $(\alpha^*, \eta^*)$, we conclude that
\beqs
\widehat V(t, \mu) =  \E\bigg[\int_t^T f(X_s^*, \alpha_s^*, \P_{X_{s}^*} ) ds  +
g(X_T^*, \P_{X_T^*})  + \int_{[t, T)} \gamma d \eta_s^* \bigg].
\enqs
\ep

\vspace{2mm}

\begin{Remark}
According to Verification Theorem \ref{veri}, if there exists $(\alpha^*, \eta^*) \in \Ac_t$ satisfying \reff{equsingularalpha*}-\reff{DPPWX*}, then $\widehat V = V$ and $(\alpha^*, \eta^*)$ is optimal. Let us denote $\Cc(V)$ and $\Dc(V)$ for the $\eta^*$-continuation region and the $\eta^*$-action region, respectively, with
\beq \label{CcV}
\Cc(V) &=& \{(t, \mu, x) \in [0, T] \times \Pc_2(\R^d) \times \R^d: \gamma + \lambda \partial_\mu V(t, \mu, x) > 0\}, \\
\Dc(V) &=& \cup_{i=1}^d \Dc_i(V),
\label{DcV}
\enq
where $$\Dc_i(V) =\{(t, \mu, x) \in [0, T] \times \Pc_2(\R^d) \times \R^d: \gamma_i + \lambda_i (\partial_\mu V)_i(t, \mu, x) = 0\}, \;\; 1 \leq i \leq d.$$
Moreover, if $(t, \P_{X_{t-}^*}, X_{t-}^*)$ starts from $\eta^*$-action region $\Dc(V)$ with positive probability, then $(t, \P_{X_{t-}^*}, X_{t-}^*)$ jumps immediately to a point on the boundary of $\Cc(V)$, denoted as $\partial \Cc(V)$. In fact, \reff{equX*nonaction} suggests that if $(t, \P_{X_{t-}^*}, X_{t-}^*)$ is in the $\eta^*$-action region with a positive probability, then $(t, \P_{X_t^*}, X_{t}^*)$ must be either in the interior of $\Cc(V)$ or in its boundary $\P$-a.s.. To check that indeed it can not be in the interior of $\Cc(V)$, let us define
\beqs
\tau: &=& \inf\Big\{ h \in [0, 1]: (t, \P_{ X_{t-}^* + h\Delta X_t^*}, X_{t-}^* + h\Delta X_t^*) \in \partial \Cc(V)\Big\}. 
\enqs
If $\P_{X_{t-}^*} \neq \P_{X_t^*}$,  according to \reff{DPPWX*}, together with the Fundamental Theorem of Calculus
\beqs
0 &=& V(t, \P_{X_t^*}) - V(t, \P_{X_{t-}^*}) + \gamma \E[\Delta \eta_t^*]\\
&=& \int_0^1 \E\Big[\Big(\lambda \partial_\mu V(t, \P_{X_{t-}^* + h\Delta X_t^*}, X_{t-}^* +  h\Delta X_t^*) + \gamma\Big) \Delta \eta_t^*\Big]dh\\
&=&  \int_0^1 \E\Big[\Big(\lambda \partial_\mu V(t, \P_{X_{t-}^* + h\Delta X_t^*}, X_{t-}^* +  h\Delta X_t^*) + \gamma\Big)1_{\{h < \tau\}} \Delta \eta_t^* \Big]d h\\
& & \;\; +  \int_0^1 \E\Big[\Big(\lambda \partial_\mu V(t, \P_{X_{t-}^* + h\Delta X_t^*}, X_{t-}^* +  h\Delta X_t^*) + \gamma\Big) 1_{\{h > \tau\}}\Delta \eta_t^*\Big]d h.
\enqs
From \reff{equX*nonaction}, the integral over $[0, \tau]$ is nonnegative  and the integral over $(\tau, 1]$ is positive, thus a contradiction unless $\tau =1$ $\P$-a.s., which means that $(t, \P_{ X_{t}^*}, X_{t}^*)$ is on $\partial \Cc(V)$ $\P$-a.s.
Similarly, if $\P_{X_t^*} = \P_{X_{t-}^*}$, according to \reff{equlinearhatVX*}, together with the Fundamental Theorem of Calculus, $\tau = 1$ $\P$-a.s.
In either case, once $(t, \P_{X_t^*}, X_t^*)$ reaches the boundary, $\eta^*_t$ acts only to prevent it from entering the interior of $\Dc(V)$.
\end{Remark}

\subsection{Relation to the Maximum Principle}

As discussed in the introduction, the maximum principle and the dynamic programming principle are the two most common approaches in solving stochastic optimal control problems. Under certain differentiability conditions, these two principles are related in the sense that  the derivatives of the value function and the solution to the adjoint equation along the optimal state process are connected. This relationship has been discussed in \cite{YZ1999}, \cite{Pham2009} for classical diffusion processes, in \cite{FOS2004}, \cite{SW2011}, \cite{ZES2012} for diffusion with jumps, and in \cite{BCM2012} for singular stochastic optimal control problems.

In this section, we will build a similar relationship for McKean--Vlasov (MKV) regular-singular control case. In order to put it in the framework suitable for the maximum principle, we will make the following assumptions.

\medskip

\noindent ({\bf H5})\; Measurable functions $b(\cdot, a, \cdot), \sigma(\cdot, a, \cdot)$ and $f(\cdot, a, \cdot)$ are differentiable with respect to $x$ and $\mu$ and all derivatives are bounded and Lipschitz continuous, where Lipschitz constants are independent of $a \in \R^m$. Moreover,  functions $b, \sigma$ and $f$ are continuously differentiable with respect to the control $a$, and all their derivatives are continuous and bounded. Similarly, $g(\cdot, \cdot)$ are differentiable with respect to $x$ and $\mu$ and all derivatives are bounded and Lipschitz continuous.

\medskip

\noindent ({\bf H6}) The Hamiltonian $H( \cdot, \cdot, \cdot, p, M)$ defined in \reff{operatorHsingular} is convex with respect to $(x, \mu, a)$, i.e., for every $x, x' \in \R^d$, $\mu, \mu' \in \Pc_2(\R^d)$, and $a, a' \in \R^m$,
\beqs
 H(x', a', \mu', p, M) &\geq& H( x, a, \mu, p, M) + \partial_x H( x, a, \mu, p, M).(x'-x)\\
  & & \;+\; \partial_a H( x, a, \mu, p, M).(a'-a) + \E\big[\partial_\mu H( x, a, \mu, p, M, \xi). (\xi' -\xi)\big],
\enqs
where $\xi$ and $\xi'$ are square integrable random variables with distributions $\mu$ and $\mu'$ (respectively).

The $\F$-adapted  processes $(p_t, M_t)_t$ are given by the adjoint equation
\beq \label{equadjoint}
\left\{
\begin{array}{rcl}
dp_t &=& -\Big\{\partial_x H(X_t^*, \alpha^*_t, \P_{X_t^*},  p_t, M_t) +\bar \E\big[\partial_\mu H(\bar X_t^*, \bar\alpha^*_t, \P_{X_t^*},  \bar p_t, \bar M_t, X_t^*)\big]\Big\}dt  \\
& & \hspace{6cm} + \;  M_t dW_t,\\
p_T &=&  \partial_x g(X_T^*, \P_{X_T^*}) + \bar\E\big[\partial_\mu g(\bar X_T^*, \P_{X_T^*}, X_T^*)\big].
\end{array}
\right.
\enq

Now recall from \cite{HSSD2018} the following maximum principle for the MKV mixed regular-singular control problem.

\begin{Proposition}

\noindent (1) Let $X^*$ be the optimal solution associated with the optimal strategy $(\alpha^*, \eta^*)$ of McKean--Vlasov mixed regular-singular control problem \reff{equdynamicssingular}-\reff{MKVsingularcost}. Assume condition ({\bf H5}). Then there exists a unique pair of $\F$-adapted  processes $(p_t, M_t)_t$ for the adjoint equation \reff{equadjoint} such that for all $(\alpha, \eta)$ $\in$ $\Ac_0 \times \Uc_0$, the following inequality holds $\P$-a.s., $dt$-a.e.
    \beqs
    \partial_a H(X_t^*, \alpha_t^*, \P_{X_t^*}, p_t, M_t)(\alpha_t - \alpha_t^*) + \E\Big[\int_{[0, T)} (\gamma + \lambda p_t) d(\eta_t - \eta^*_t)\Big] \geq 0.
    \enqs

\noindent (2) Let $(\hat\alpha, \hat\eta) \in \Ac_0 \times \Uc_0$ be an admissible control, and $X_t^{\hat\alpha, \hat\eta}$, $(p^{\hat \alpha}_t, M^{\hat\alpha}_t)$ be the solution of \reff{equdynamicssingular} and \reff{equadjoint} respectively associated with $(\hat\alpha, \hat\eta)$. Assume conditions ({\bf H5})-({\bf H6}). And suppose that the singular control $(\hat\alpha, \hat\eta)$ satisfies that
    \beqs
    \E \Big[\int_0^T \partial_a H(X_t^{\hat\alpha, \hat\eta}, \hat\alpha_t, \P_{X_t^{\hat\alpha, \hat\eta}}, p_t^{\hat\alpha}, M_t^{\hat\alpha})(\alpha_t - \hat\alpha_t)dt\Big] + \E\Big[\int_{[0, T)} (\gamma + \lambda p_t^{\hat\alpha}) d(\eta_t - \hat\eta_t)\Big] \geq 0,
    \enqs
    for any $(\alpha, \eta) \in \Ac_0 \times \Uc_0$.
    Then $(\hat\alpha, \hat\eta)$ is an optimal control of the McKean--Vlasov singular control problem \reff{equdynamicssingular}-\reff{MKVsingularcost} starting at time $t$ $=$ $0$.
\end{Proposition}

We now show that the adjoint process can be expressed in terms of the Lions derivative of the value
function for the MKV mixed regular-singular control problem.

\begin{Theorem} \label{Thmrelationadjoint}
Assume that $V \in \Cc^{1, (2,1)}([0, T] \times \Pc_2(\R^d))$.
If $(\alpha^*, \eta^*) \in \Ac_0 \times \Uc_0$ satisfies conditions \reff{equX*nonaction}-\reff{DPPWX*}, then the solution of the adjoint equation \reff{equadjoint} is given by
\beq \label{relationpq}
p_t = \partial_\mu V(t, \P_{X_t^*}, X_t^*),\;\; M_t = \partial_{x}\partial_\mu V(t, \P_{X_t^*}, X_t^*)\sigma( X_t^*, \alpha^*_t, \P_{X_t^*}),  \quad 0\leq t\leq T.
\enq
\end{Theorem}
{\bf Proof.}\;
Observe that
\beqs
& & \min_{(\mu, a) \in \Pc_2(\R^d) \times \R^m} F(t, \mu, a)= F(t, \P_{X_t^*}, \hat a(t, X_t^*, \P_{X_t^*})),
\enqs
where
$F(t, \mu, a)$ is given by
\beqs
F(t, \mu, a): &=&  \partial_t V(t, \mu) + \E\big[H(\xi, a, \mu, \partial_\mu V(t, \mu, \xi), \partial_x\partial_\mu V(t, \mu, \xi))\big].
\enqs
Differentiating $F(t, \mu, a)$ with respect to $\mu$ and then evaluating at $\P_{X_t^*}$ and $\hat a$ yields
\beq \label{equpartialmu}
\partial_\mu F(t, \P_{X_t^*}, \hat a(t, X_t^*, \P_{X_t^*})) =0.
\enq
To simplify notations, set
\beqs
b_t^*: = b(X_t^*, \alpha^*_t, \P_{X_t^*}), \;\bar b_t^*: = b(\bar X_t^*, \bar \alpha^*_t, \P_{X_t^*}), \;\sigma_t^*: = \sigma(X_t^*, \alpha^*_t, \P_{X_t^*}), \; \bar\sigma_t^*: =\sigma(\bar X_t^*, \bar \alpha^*_t, \P_{X_t^*}).
\enqs
Then straightforward calculation of \reff{equpartialmu} suggests
\beq
& & \partial_t\partial_\mu V(t, \P_{X_t^*}, X_t^*)  +  \partial_x\partial_\mu V(t, \P_{X_t^*}, X_t^*)b_t^*+ \bar \E\big[ \partial_\mu^2 V(t, \P_{X_t^*}, \bar X_t^*, X_t^*) \bar b_t^*\big]\nonumber\\
& & \; + \; \frac{1}{2} \Big \lbrace {\rm Tr}\big(\sigma_t^*(\sigma_t^*)\trans\partial_{xx}\partial_\mu V(t, \P_{X_t^*}, X_t^*)\big) + \bar \E\big[ {\rm Tr}\big(\bar\sigma_t^*(\bar\sigma_t^*)\trans \partial_x\partial_\mu^2 V(t, \P_{X_t^*}, \bar X_t^*, X_t^*)\big)\big]\Big \rbrace \nonumber\\
&=& - \Big\lbrace \partial_x H(X_t^*, \alpha_t^*, \P_{X_t^*}, \partial_\mu V(t, \P_{X_t^*}, X_t^*), \partial_x\partial_\mu V(t, \P_{X_t^*}, X_t^*)\sigma_t^*) \nonumber\\
& & \; + \;\bar \E\big[ \partial_\mu H(\bar X_t^*, \bar\alpha_t^*, \P_{X_t^*}, \partial_\mu V(t, \P_{X_t^*}, \bar X_t^*), \partial_x\partial_\mu V(t, \P_{X_t^*}, \bar X_t^*), X_t^*)\bar\sigma_t^*)\big] \Big \rbrace.
\label{equderivativeMP}
\enq
Applying Corollary \ref{corollary2} to $\partial_\mu V(s, \P_{X_s^*}, X_s^*)$ between $t$ and $T$ yields
\beq \label{ItopartialmuV}
& & \partial_\mu V(t, \P_{X_t^*}, X_t^*) - \partial_\mu V(T, \P_{X_{T}^*}, X_{T}^*) \\
& =&  - \int_t^T \biggl\lbrace \partial_s\partial_\mu V(s, \P_{X_s^*}, X_s^*) + \partial_x\partial_\mu V(s, \P_{X_s^*}, X_s^*)b_s^*  + \frac{1}{2} {\rm Tr} \big(\sigma^*_s(\sigma^*_s)\trans\partial_{xx}\partial_\mu V(s, \P_{X_s^*}, X_s^*) \big)\biggl\rbrace ds \nonumber\\
& &-  \int_t^T \bar \E \Big[\partial_{\mu}^2 V(s, \P_{X_s^*}, X_s^*, \bar X_s^*) \bar b_s^* +  \frac{1}{2} {\rm Tr} \big(\bar\sigma_s^*(\bar\sigma_s^*)\trans\partial_{\bar x} \partial_{\mu}^2V(s, \P_{X_s^*}, X_s^*, \bar X_s^*) \big)\Big] ds \nonumber\\
 & & - \int_t^T \partial_x \partial_\mu V(s, \P_{X_s^*}, X_s^*) \sigma_s^* dW_s \nonumber\\
 & & - \lambda \underbrace{\Big(\int_t^T \partial_x \partial_\mu V(s, \P_{X_s^*}, X_s^*)d\eta^*_s - \sum_{t \leq s < T}\partial_x \partial_\mu V(s, \P_{X_s^*}, X_s^*)\Delta \eta^*_s \Big)}_{I_1}\nonumber\\
& &  - \underbrace{\sum_{t \leq  s < T} \Big(\partial_\mu V(s, \P_{X_s^*}, X_s^*) - \partial_\mu V(s, \P_{X_s^*}, X_{s-}^*) \Big)1_{\{\P_{X_{s}^*} = \P_{X_{s-}^*}\}}}_{I_{2}^1}  \nonumber\\
& & -  \underbrace{\sum_{t \leq  s < T} \Big(\partial_\mu V(s, \P_{X_s^*}, X_s^*) - \partial_\mu V(s, \P_{X_s^*}, X_{s-}^*) \Big)1_{\{\P_{X_{s}^*} \neq \P_{X_{s-}^*}\}}}_{I_{2}^2}\nonumber\\
& & -\lambda \underbrace{\bar \E \Big[\int_t^T \partial_\mu^2  V(s, \P_{X_s^*}, X_s^*, \bar X_s^*) d\bar \eta^*_s - \sum_{t \leq s < T} \partial_\mu^2 V(s, \P_{X_s^*}, X_s^*, \bar X_s^*)\Delta \bar\eta^*_s\Big]}_{I_3} \nonumber\\
& & -    \underbrace{\bar\E\Big[ \sum_{t \leq s < T} \Big(\frac{\delta (\partial_\mu V)}{\delta \mu}(s, \P_{X_s^*, }, X_s^*, \bar X_{s}^*) -  \frac{\delta (\partial_\mu V)}{\delta \mu}(s, \P_{X_{s}^*}, X_s^*, \bar X_{s-}^*) \Big) 1_{\{\P_{X_s^*} = \P_{X_{s-}^*}\}}\Big]}_{I_4}  \nonumber\\
& &  -  \underbrace{\sum_{t \leq s < T} \Big(\partial_\mu V(s, \P_{X_s^*}, X_{s-}^*) - \partial_\mu V(s, \P_{X_{s-}^*}, X_{s-}^*)\Big)1_{\{\P_{X_{s}^*} \neq \P_{X_{s-}^*}\}}}_{I_5} .\nonumber
\enq

Now, we calculate $I_1$ and $I_3$. On one hand, $\Dc(V) \subset \arg\min\{\gamma + \lambda \partial_\mu V(\mu, x)\}$, then for every $(t, \mu, x) \in \Dc(V)$ in \reff{DcV},
\beqs
\partial_\mu \Big(\gamma + \lambda \partial_\mu V(\mu, x)\Big) =\partial_{\mu\mu}V(\mu, x, x') = 0,\;\; \partial_x \Big(\gamma + \lambda \partial_\mu V(\mu, x)\Big) =\partial_x\partial_\mu V(\mu, x) =0.
\enqs
On the other hand, \reff{equeta*} implies that $\E\Big[\int_t^T 1_{(s, \P_{X_s^*}, X_s^*) \in \Cc(V)} d \eta_s^{*, c}\Big] =0$, therefore
\beqs
I_1: &=& \lambda\int_t^T \partial_x\partial_\mu  V(s, \P_{X_s^*}, X_s^*) d\eta^{*, c}_s\\
&=& \lambda\int_t^T  1_{\{(s, \P_{X_s^*}, X_s^*) \in \Cc(V)\}}\partial_x\partial_\mu  V(s, \P_{X_s^*}, X_s^*) d \eta^{*, c}_s\\
& &\;+ \; \lambda\int_t^T 1_{\{(s, \P_{X_s^*}, X_s^*) \in \Dc(V)\}}\partial_x\partial_\mu  V(s, \P_{X_s^*}, X_s^*)d \eta^{*, c}_s \\
&=& 0.
\enqs
Similarly, we have
\beqs
I_3: =\bar \E \Big[\int_t^T \partial_\mu^2  V(s, \P_{X_s^*}, X_s^*, \bar X_s^*) d\bar \eta^*_s - \sum_{t \leq s < T} \partial_\mu^2 V(s, \P_{X_s^*}, X_s^*, \bar X_s^*)\Delta \bar\eta^*_s\Big] =0.
\enqs
Next, we calculate $I_2^2 + I_5$. By the Fundamental Theorem of Calculus,
\beqs
I_2^2 + I_5 : &=& \sum_{t \leq s < T}\Big(\partial_\mu V(s, \P_{X_s^*}, X_s^*) - \partial_\mu V(s, \P_{X_{s-}^*}, X_{s-}^*)\Big) 1_{\{\P_{X_{s}^*} \neq \P_{X_{s-}^*}\}}\\
 &=& \lambda \sum_{t \leq s < T} \Big(\int_0^1 \partial_x\partial_\mu V(s, \P_{X_{s-}^* + h \Delta X_s^*}, X_{s-}^* + h \Delta X_s^*) \Delta \eta_s^* d h \\
& +&  \int_0^1 \bar\E\big[\partial_\mu^2 V(s, \P_{X_{s-}^* + h \Delta X_s^*}, X_{s-}^* + h \Delta X_s^*, \bar X_{s-}^* + h \Delta \bar X_r^*) \Delta \bar \eta_s^*\big] d h\Big) 1_{\{\P_{X_{s}^*} \neq \P_{X_{s-}^*}\}}.
\enqs
To prove that the right hand side in the above equation vanishes, it suffices to check that if $\Delta \eta_s^* > 0$, then for every $h \in [0, 1]$, $\P$-a.s., $\bar \P$-a.s.
\beq \label{equzeropartialxmuV}
& & \partial_x\partial_\mu V(s, \P_{X_{s-}^* + h \Delta X_s^*}, X_{s-}^* + h \Delta X_s^*) = 0,\\ & & \partial_\mu^2 V(s, \P_{X_{s-}^* + h \Delta X_s^*}, X_{s-}^* + h \Delta X_s^*, \bar X_{s-}^* + h \Delta \bar X_s^*) =0.
\label{equzeropartialmu2V}
\enq
From \reff{DPPWX*}, we see that when $s \in$ $\{t \leq s < T: \P_{X_s^*} \neq \P_{X_{s-}^*}\}$,
\beqs
0 &=&  V(s, \P_{X_s^*}) - V(s, \P_{X_{s-}^*})  + \gamma \mathbb{E}  \big[\Delta \eta_s^* \big]\\
&=& \int_0^1 \E\Big[\Big(\lambda \partial_\mu V(s, \P_{X_{s-}^* + h \Delta X_s^*}, X_{s-}^* + h \Delta X_s^*) + \gamma \Big) \Delta \eta_s^*\Big]dh.
\enqs
The right hand side is nonnegative, $\Delta \eta_s^* > 0$ implies that $(s, \P_{X_{s-}^* + h \Delta X_s^*}, X_{s-}^* + h \Delta X_s^*) \in \Dc(V)$, $\P$-a.s., for every $ h \in [0, 1]$. This implies \reff{equzeropartialxmuV}-\reff{equzeropartialmu2V}.
Hence, $I_2^2 + I_5 = 0$.

Now we compute $I_2^1 + I_4$. By the Fundamental Theorem of Calculus,
\beqs
I_2^1 &=& \sum_{t \leq  s < T} \Big(\int_0^1 \partial_x\partial_\mu V(s, \P_{X_{s}^*}, X_{s-}^* + h \Delta X_s^*) \Delta \eta_s^* dh\Big) 1_{\{\P_{X_{s}^*} = \P_{X_{s-}^*}\}}, \\
I_4  &=& \bar\E\Big[\sum_{t \leq s < T} \Big(\int_0^1 \lambda \partial_{\mu}^2V (s, \P_{X_{s-}^* + h \Delta X_s^*}, X_s^*, \bar X_{s-}^* + h \Delta \bar X_s^*) \Delta \eta_s^* dh\Big)1_{\{\P_{X_s^*} = \P_{X_{s-}^*}\}}\Big].
\enqs
To prove the right hand side in the above equation vanishes, it suffices to show if $\Delta\eta_s^* > 0$, then for every $h \in [0, 1]$, $\P$-a.s.,  $\bar\P$-a.s.
\beq \label{equzeropartialxmu2V}
\partial_x\partial_\mu V(s, \P_{X_s^*}, X_{s-}^* + h \Delta X_s^*) = 0, \;\; \partial_{\mu}^2V (s, \P_{X_s^*}, X_s^*, \bar X_{s-}^* + h \Delta \bar X_s^*) = 0.
\enq
From \reff{equlinearhatVX*},  when $s \in$ $\{t \leq s < T: \P_{X_s^*} = \P_{X_{s-}^*}\}$
\beqs
0 &=& \frac{\delta V}{\delta \mu}(s, \P_{X_{s}^*}, X_{s}^*) - \frac{\delta V}{\delta \mu} (s, \P_{X_{s}^*}, X_{s-}^*)  + \gamma \Delta \eta_s^* \\
&=& \int_0^1 \Big(\lambda \partial_\mu V(s, \P_{X_s^*}, X_{s-}^* + h \Delta X_s^*)  + \gamma \Big) \Delta \eta_s^* dh.
\enqs
The right hand side is nonnegative, $\Delta \eta_s^* > 0$ implies that $(s, \P_{X_{s}^*}, X_{s-}^* + h \Delta X_s^*) \in \Dc(V)$, $\P$-a.s., for every $ h \in [0, 1]$. This implies \reff{equzeropartialxmu2V}
\beqs
\partial_x\partial_\mu  V(s, \P_{X_s^*}, X_{s-}^* + h \Delta X_s^*)  = 0, \;\;\; \partial_\mu^2 V(s, \P_{X_s^*}, X_{s-}^* + h \Delta X_s^*, \bar X_s^*) = 0.
\enqs
Hence $I_2^1 =0$ and $I_4=0$.

Note that all terms $I_1-I_5$ associated with $\eta^*_s$ vanish. By \reff{equderivativeMP}-\reff{ItopartialmuV}, we obtain
\beqs
\left\{
\begin{array}{rcl}
d \partial_\mu V(t, \P_{X_t^*}, X_t^*) &=& - \Big\lbrace \partial_x H(X_t^*, \alpha_t^*, \P_{X_t^*}, \partial_\mu V(t, \P_{X_t^*}, X_t^*), \partial_x\partial_\mu V(t, \P_{X_t^*}, X_t^*)\sigma_t^* ) \nonumber\\
& + & \bar \E\big[ \partial_\mu H(\bar X_t^*, \bar\alpha_t^*, \P_{X_t^*}, \partial_\mu V(t, \P_{X_t^*}, \bar X_t^*), \partial_x\partial_\mu V(t, \P_{X_t^*}, \bar X_t^*)\bar\sigma_t^*)\big] \Big \rbrace dt\\
& + &  \partial_x\partial_\mu V(t, \P_{X_t^*}, X_t^*) \sigma(X_t^*, \alpha_t^*, \P_{X_t^*}) dW_t,\\
\partial_\mu V(T, \P_{X_T^*}, X_T^*) &=& \partial_x g(X_t^*, \P_{X_t^*}) + \bar\E[\partial_\mu g(\bar X_T^*, \P_{X_T^*}, X_T^*)].
\end{array}
\right.
\enqs
Thus, the pair $(p, M)$ given by \reff{relationpq} is the solution of the adjoint equation \reff{equadjoint}.
\ep

\begin{Remark}[Mean-field games with singular controls]
\label{remMFGMKV} When the controls are regular,  mean-field game (MFG) and McKean--Vlasov control are related (\cite{CD2018I}, chapter 6).  In the case of singular controls, \cite{FH2017} established the existence of an optimal control to McKean--Vlasov singular control using a similar method for MFGs with singular controls. In fact, one may mimic the case of regular controls and  show that MFG with singular control and McKean--Vlasov singular control are also connected. To see this, assume for simplicity that $\sigma$ in \reff{equdynamicssingular} does not depend upon on the regular control $\alpha$, then the FBSDE system given by the maximum principle \reff{equdynamicssingular}-\reff{equadjoint} for the McKean--Vlasov control problem may be also identified with (at least formally) the FBSDE system given by applying the maximum principle to the following auxiliary MFG problem, where the dynamics follow
\begin{align*}\label{equdynamicssingularMFG}
\left\{
\begin{array}{rcl}
d{X_s}  &= &{b} (X_s, \alpha_s, \mu_s) ds + {\sigma}(X_s, \mu_s)  d{W_s} +  {\lambda} d\eta_s, \; t \leq s \leq T,\\
 X_{t-} &=& \xi \in L^2(\Fc_t; \R^d),
 \end{array}
 \right.
\end{align*}
and the cost functional is given by
\beqs
I(\alpha, \eta) &=& \E\biggl[g(X_T, \mu_T) + \int_{\R^d} \frac{\delta g}{\delta \mu}(y, \mu_T, X_T)\mu_T(dy) + \int_t^T f(X_s, \mu_s, \alpha_s)ds\\
&  + & \int_t^T \int_{\R^d} \frac{\delta H}{\delta \mu}\Big(y, \hat\alpha\big(\mu_s, \partial_y u(s, y), \partial_{yy}^2 u(s, y)\sigma(y, \mu_s)\big), \mu_s, \partial_y u(s, y), \partial_{yy}^2 u(s, y)\sigma(y, \mu_s)\Big) \mu_s(dy) \\
& + & \int_t^T \gamma d \eta_s \biggl].
\enqs
Here $\hat\alpha(x, \mu, p, M)$ is the minimizer of $H(x, a, \mu, p, M)$ in \reff{operatorHsingular}, and the function $u$ is the linear derivative of the value function $V(t, \mu)$ in \reff{equDPP} along the optimal path.
\end{Remark}

\subsection{Example: Mean-Variance Singular Control}

We now analyze a class of  one-dimensional mean-variance singular control problem, where coefficients of the dynamics in \reff{equdynamicssingular} are specialized with
\beqs
b (x, a, \mu) &=& rx + \rho a, \\
\sigma(x, a, \mu) &=& \sigma a,
\enqs
for $(x, \mu, a) \in \R \times \Pc_2(\R) \times \R$, with $r, \rho$ and $\sigma >0$ constants in $\R$. And in \reff{MKVsingularcost}, the running cost $f$ $\equiv$ $0$,
and the mean-variance terminal cost function is
\beqs
g(x, \mu) &=& \frac{\beta}{2}(x - \bar \mu)^2- x,
\enqs
where the constant $\beta >0$, and
$\bar\mu: = \int_{\R} x \mu(dx)$.

We will search for a classical solution to the dynamic programming equation \reff{equsingularDP}. In the waiting region $\Cc(V)$ in \reff{CcV}, $V(t, \mu)$ would satisfy the HJB equation \reff{equHJBnonaction}, which corresponds to the classical linear quadratic McKean--Vlasov control problem.  Now take $V(t, \mu)$ of the following form
\beqs
V(t, \mu) = A(t){\rm Var}(\mu) + B(t)\bar \mu^2 + C(t) \bar\mu + D(t),
\enqs
 for some time-dependent functions $A$, $B$, $C$, and $D: [0, T] \to \R$.  Solving the corresponding HJB equation \reff{equHJBnonaction} in the similar way as in Section \ref{secMKVLQ}, we see that $A(t), B(t), C(t)$ and $D(t)$ satisfy the following ODEs
\beqs
\left\{
\begin{array}{rcl}
\dot A(t) - \Big(\frac{\rho^2}{\sigma^2} - 2 r\Big) A(t) &=&0,\;\;\; A(T) =\frac{\beta}{2},\\
\dot B(t) - \frac{\rho^2}{\sigma^2} \frac{B^2(t)}{A(t)} +  2 r B(t) &=&0,\;\;\; B(T) = 0,\\
\dot C(t) + r C(t) -\frac{\rho^2}{\sigma^2} \frac{B(t)}{A(t)} &=&0,\;\;\; C(T) =-1,\\
\dot D(t) - \frac{\rho^2}{\sigma^2} \frac{C(t)^2}{4 A(t)} &=&0,\;\;\; D(T) =0,
\end{array}
\right.
\enqs
which can be explicitly solved such that
\beqs
A(t) &=& \frac{\beta}{2}\exp \Big((2 r - \frac{\rho^2}{\sigma^2}) (T -t)\Big),\; B(t)=0,\\
C(t)&=&- \exp\Big(r(T-t)\Big),\; D(t) =\frac{1}{4} \exp\Big(\frac{\rho^2}{\sigma^2} (T-t)\Big) -\frac{1}{4}.
\enqs
Therefore, $\Cc(V)$ in \reff{CcV} and $\Dc(V)$ in \reff{DcV} are now given by
\beqs
\Cc(V) &=& \Big\{(t, \bar \mu, x) \in [0, T] \times \R \times \R : {\lambda\beta}\exp \big((2 r - \frac{\rho^2}{\sigma^2}) (T -t)\big) (x- \bar\mu)\\
& &  \hspace{5cm} -\lambda \exp\big(-r (T-t)\big) + \gamma > 0\Big\},\\
\Dc(V) &=& \Big\{(t, \bar \mu, x) \in [0, T] \times \R \times \R : {\lambda\beta}\exp \big((2 r - \frac{\rho^2}{\sigma^2}) (T -t)\big) (x- \bar\mu)\\
& &  \hspace{5cm} -\lambda \exp\big(-r (T-t)\big) + \gamma = 0\Big\}.
\enqs
From \reff{equX*nonaction}-\reff{DPPWX*},
$\P( \text {Leb a.e. } t \in [0, T], \; (t, \P_{X_t^*}, X_t^*) \in \Cc(V)) =1,$
where $X_t^*$ is the controlled process associated with the optimal control $(\alpha^*, \eta^*)$, where
\beqs
\alpha^*_t = -\frac{\rho}{\sigma^2} (X_t^* - \E[X_t^*]) + \frac{\rho}{\beta\sigma^2} \exp \Big(( \frac{\rho^2}{\sigma^2}- r) (T -t)\Big),
\enqs
and $\eta^*$
\beq \label{MVexampleeta*}
\E \int_0^T 1_{\{(t, \P_{X_t^*}, X_t^*) \in \Dc(V)\}} d \eta^*_t =0.
\enq
The controlled process $X_t^*$ then follows
\beq \label{MVexampleX*}
dX_t^* &=& \Big((r -\frac{\rho^2}{\sigma^2}) X_t^* + \frac{\sigma^2}{\rho^2} \E[X_t^*] + \frac{\rho^2}{\beta\sigma^2} \exp(\frac{\rho^2}{\sigma^2} - r)(T-t))\Big)dt\\
& &\; +\;\Big(-\frac{\rho}{\sigma} (X_t^* -\E[X_t^*]) + \frac{\rho}{\beta\sigma}\exp(\frac{\rho^2}{\sigma^2} - r)(T-t)) \Big) dW_t + \lambda d\eta^*_t. \nonumber
\enq
Solving  \reff{MVexampleeta*}-\reff{MVexampleX*} is equivalent to solving  a one-dimensional Skorokhod problem,  see \cite{HOS2017}. By Theorem 4.1 in \cite{HOS2017}, such an $\eta^*$ exists uniquely.





\bibliographystyle{plain}

\bibliography{refs}

\end{document}